\documentclass[a4paper]
{cedram-smai-jcm}

\usepackage{tikz}
\usepackage{subcaption}
\usepackage{graphicx} 
\usepackage{geometry}
\usepackage{amsmath,amssymb,amsthm,mathtools}
\usepackage{mathrsfs}
\usepackage{subcaption}
\usepackage{hyperref}
\usepackage{todonotes}
\usepackage{multirow}
\usepackage{algorithm}
\usepackage{algpseudocode}
\usepackage{booktabs}   
\usepackage{multirow}   
\usepackage{rotating}
\usepackage{lscape}

\newcommand{\R}{\mathbb{R}}

\title[Low-rank MLMC for Kinetic Equations]{A Dynamical Low-rank Multilevel Monte Carlo Estimator for High-Dimensional Kinetic Equations}

\author[C. Patwardhan]{\firstname{Chinmay} \lastname{Patwardhan}}
\address{Dept. of Mathematics, Karlsruhe Institute of Technology, Karlsruhe, Germany}
\email{chinmay.patwardhan@kit.edu}
\thanks{Chinmay Patwardhan was funded by the Deutsche Forschungsgemeinschaft (DFG, German Research Foundation) – Project-ID 258734477 – SFB 1173}

\author[S. Krumscheid]{\firstname{Sebastian} \lastname{Krumscheid}}
\address{Scientific Computing Center, Karlsruhe Institute of Technology, Karlsruhe, Germany}
\email{sebastian.krumscheid@kit.edu}

\author[J. Kusch]{\firstname{Jonas} \lastname{Kusch}}
\address{Scientific Computing, Norwegian University of Life Sciences, \AA s, Norway}
\email{jonas.kusch@nmbu.no}

\author[E. L\o vbak]{\firstname{Emil} \lastname{L\o vbak}}
\address{Scientific Computing Center, Karlsruhe Institute of Technology, Karlsruhe, Germany}
\email{emil.loevbak@kit.edu}
\thanks{Emil Løvbak was funded by the Deutsche Forschungsgemeinschaft
(DFG, German Research Foundation) – Project-ID 563450842}

\author[P. Stammer]{\firstname{Pia} \lastname{Stammer}}
\address{Dept. of Radiation Science \& Technology, Delft University of Technology (TU Delft), 
Delft, Netherlands}
\email{p.k.stammer@tudelft.nl}
\thanks{Pia Stammer received funding from the German National Academy of Sciences Leopoldina for the project underlying this article, under grant number LPDS 2024-03}

  

\begin{document}

\begin{abstract}
Kinetic equations are used to model a wide range of phenomena important for real-world applications. Their applications span astrophysics, nuclear physics, engineering, and social sciences. Due to their high-dimensional phase space, modelling and quantifying uncertainties, relevant for applications, poses a significant challenge even for modern computing infrastructure. In recent years, dynamical low-rank approximation (DLRA) has gained popularity for making fine grid simulations of high-dimensional problems feasible by evolving the solution of a time-dependent PDE as a low-rank factorization. This reduces the computational and memory requirements significantly.

In this work, we propose a low-rank multilevel Monte Carlo estimator for kinetic equations based on a probabilistic rank-adaptive DLRA time integrator. The level hierarchy of the low-rank multilevel estimator is constructed through spatial refinement and by ensuring that the low-rank error remains below the spatial discretization error. We demonstrate the efficacy of the estimator through several numerical experiments from radiation transport, radiation therapy, and shallow water flow.
\end{abstract}

\maketitle


\section{Introduction}

Kinetic equations model the microscopic effects of particle transport and their interactions, such as collisions and absorption, through a probability density function. By bridging the gap between atomistic particle models and macroscopic continuum models, kinetic equations play a crucial role in modelling multi-scale dynamics of transport phenomena. These equations describe a wide range of physical phenomena ranging from gas dynamics and radiation transport to shallow water flow and chemical reaction networks, with applications spanning astrophysics, nuclear physics, radiation therapy, engineering, and social sciences~\cite{prugger_dynamical_2023,lewis_computational_1984,armbruster_kinetic_2003,kowalski_moment_2019}. Since real-world applications often introduce uncertainties, including these in models provides a more realistic description. These uncertainties stem from modeling errors, including the use of empirical collision kernels; measurement errors, such as incorrect measurements of the physical domain; and device errors, including errors in initial and boundary data. Quantifying these uncertainties is thus essential for robustness and reliability in real-world applications.

In this work, we consider uncertain kinetic equations of the type 
\begin{subequations}\label{eq:genericKEwU}
    \begin{align}
        \begin{split}
            \partial_{t}u(t,x,v\,;\omega) + \mathcal{A}(u)\partial_{x}u(t,x,v\,;\omega)  &= \mathcal{S}(u\,;\omega)(t,x,v\,;\omega)
        \end{split}\label{eq:genericKEwUMain}\\
        \begin{split}
            u(0,x,v\,;\omega) &= u_{0}(x,v\,;\omega),
        \end{split}\label{eq:genericKEwUIC}\\
        \begin{split}
         u(t,x,v\,;\omega) &= u_{b}(t,x,v\,;\omega),\quad x\in\Gamma_{\mathcal{D}},   
        \end{split}\label{eq:genericKEwUBC}
    \end{align}    
\end{subequations}
where $u(t,x,v\,;\omega)$ denotes the phase space density at time $t\in\mathbb{R}_{\geq 0}$, position $x=(x_{1},\ldots,x_{d})\in\mathcal{D}\subseteq\mathbb{R}^{d}$, moving with velocity $v=(v_{1},\ldots,v_{d})\in\mathcal{V}\subseteq\mathbb{R}^{d}$. The random variable $\omega\in\mathbb{R}^{k}$ models the uncertainty and has the probability density function (pdf) $p(\omega)$. The left-hand side of \eqref{eq:genericKEwUMain}, known as the streaming term, models particle movement in the phase-space where $\mathcal{A}(u)\partial_{x}u$ is known as the transport or advection term. Changes in the velocity of particles due to inter-particle collisions or scattering and absorption by the background medium are modelled by the collision operator $\mathcal{S}$. The initial condition is given by $u_{0}$ and the boundary conditions are given by $u_{b}$ on the boundary $\Gamma_{\mathcal{D}}$ of $\mathcal{D}$. The streaming and scattering operators are problem-dependent and determine the dynamics of the system along with the initial and boundary conditions. 

We are often interested in quantifying uncertainties in the entire solution of a PDE or a quantity derived from it. Let $Q : \omega\mapsto Q(\omega)=Q(u(\cdot,\cdot,\cdot\,;\omega))\in H $ be a bounded linear functional or a Lipschitz function of the weak solution $u(t,x,v\,;\omega)$ mapping into the Hilbert space $H$. Then $Q$ denotes the quantity of interest (QoI) of the solution. For instance, if $Q$ is a functional of the solution, then $H = \mathbb{R}$, while if it is the weak solution at some time $T\geq0$, then $Q(\omega) := u(T,\cdot\,,\cdot\,;\omega)\in H= L^{2}(\mathcal{D}\times\mathcal{V})$. In this work, we are interested in estimating the expected value of the QoI under the random variable $\omega\sim p(\omega)$, i.e.
\begin{equation}\label{eq:ExpectedVal}
    \mathbb{E}\left[Q\right] \coloneqq \int_{\mathbb{R}^{k}} Q(u(\cdot,\cdot,\cdot\,;\omega))\,p(\omega)\mathrm{d}\omega, 
\end{equation}
with the variance
\begin{equation}\label{eq:VarianceVal}
    \mathbb{V}[Q] \coloneqq \mathbb{E}\left[\left\lVert Q(\omega) - \mathbb{E}\left[Q\right]\right\lVert_{H}^{2}\right].
\end{equation}

There are several ways of estimating $\mathbb{E}[Q]\in H$, broadly classified as \textit{intrusive} and \textit{non-intrusive} methods. As the name suggests, intrusive uncertainty quantification (UQ) methods, such as the generalized polynomial chaos (gPC)-based stochastic Galerkin (sG) method, modify the underlying equations to obtain estimates of $\mathbb{E}\left[Q\right]$. In contrast, non-intrusive methods, such as sampling-based Monte Carlo (MC) methods and collocation-based methods, do so without altering the governing equations. A major challenge to quantifying uncertainties in kinetic equations is the high computational complexity and memory requirements for grid-based simulations. This is due to its high-dimensional phase space, which for each realization of $\omega$, including time, is at least $2d+1$ dimensional, where $d$ is the spatial dimension. In the literature, this phenomenon is referred to as the \textit{curse of dimensionality}, wherein computational costs grow exponentially with the number of dimensions. Mitigation techniques, such as reduced-order models~\cite{liu_bi-fidelity_2020,pareschi_introduction_2021}, surrogate models~\cite{chen_structure_2025,jin_asymptotic-preserving_2024}, and low-rank approximations~\cite{fairbanks_low-rank_2017,koch_dynamical_2007,einkemmer_review_2025,kusch_dynamical_2022,patwardhan_low-rank_2026}, have been proposed to reduce the computational complexity for probabilistic and deterministic kinetic equations.  

In recent years, dynamical low-rank approximation (DLRA)~\cite{koch_dynamical_2007} has gained significant popularity in the kinetic theory community as an efficient numerical method for approximating deterministic high-dimensional time-dependent PDEs. DLRA is a hybrid approach in which time-dependent low-rank factors of the solution are updated according to the problem dynamics. This approach simultaneously reduces the computational costs and memory requirements of a simulation. A recent review paper highlights the use of DLRA in the context of kinetic equations~\cite{einkemmer_review_2025}. We refer the reader to this work for a detailed comparison of the methods and \cite{einkemmer_low-rank_2018,einkemmer_asymptotic-preserving_2024,kusch_dynamical_2022,patwardhan_parallel_2026,kusch_robust_2023} for further reading on DLRA in kinetic theory.

We are interested in developing DLRA-based UQ methods for kinetic equations. The first such approach using a slightly different representation was developed under the term dynamically orthogonal (DO) field expansion in \cite{sapsis_dynamically_2009} and applied to problems in ocean dynamics with random data \cite{ueckermann2013numerical}. Apart from the DO-type methods, there are primarily two approaches in the literature that address the combination of DLRA and uncertainties. The first approach treats the random variable $\omega$ as a phase-space variable and evolves an uncertain basis over time and has been explored in \cite{musharbash_dual_2018,musharbash_symplectic_2020,kusch_dynamical_2022}. This approach can be viewed as a stochastic Galerkin method with a time-dependent, adaptive random basis. This approach can be viewed as a stochastic Galerkin method with a time-dependent, adaptive random basis. The second approach, outlined in \cite{patwardhan_low-rank_2026}, treats the random variable as an input parameter to the PDE and uses DLRA as a black-box approximation method for MC-based sampling. MC-based methods are widely adopted due to their ease of implementation and widespread use for propagating uncertainties in complex phenomena. Variance-reducing enhancements, such as the multilevel Monte Carlo (MLMC)~\cite{giles_multilevel_2008,giles_multilevel_2015} method, have further accelerated the use of sampling-based methods. MLMC estimators exploit the hierarchical structures of numerical methods to construct estimators that balance stochastic and numerical approximation errors while minimizing computational costs. Extensive work on MLMC theory and estimators for elliptic equations, parabolic equations, and hyperbolic conservation laws can be found in~\cite{cliffe_multilevel_2011,collier_continuation_2015,grote_uncertainty_2022,barth_multilevel_2013,luo_multilevel_2019,mishra_sparse_2012,mishra_multi-level_2012}. Hence, we are interested in the second approach of combining DLRA with a sampling-based UQ methods, wherein we combine non-intrusive methods with DLRA. 

To mitigate the curse of dimensionality, in our previous work \cite{patwardhan_low-rank_2026}, we constructed a DLRA-based low-rank Monte Carlo (DLR-MC) estimator for radiation transport in one-dimensional slab geometry. This estimator computed the low-rank approximation based on a \textit{fixed-rank integrator}, i.e., one rank was chosen apriori and stayed constant over time. The variance of the estimator was further reduced by constructing a control variates estimator (DLR-CV) based on a ``high" rank - ``low" rank multi-fidelity strategy. The DLR-MC and DLR-CV estimators focus on reducing the statistical error of the QoI, with the implicit assumption that it dominates the numerical error. In higher dimensions $(d\geq 2)$, this imposes a severe restriction on the extent of variance reduction since a highly resolved grid is required to satisfy this assumption. Moreover, these estimators utilize fixed-rank integrators, which necessitate prior knowledge of the rank (or maximal rank) of the solution over the entire simulation. This leads to an over-approximation of rank, resulting in higher computational costs.

In this work, we propose a rank-adaptive low-rank MC estimator based on the augmented BUG integrator~\cite{ceruti_rank-adaptive_2022}, which improves on the fixed-rank augmented BUG estimator used in \cite{patwardhan_low-rank_2026}. Further, we design a novel rank-adaptive dynamical low-rank MLMC framework for kinetic equations. The central challenge is efficiently balancing the low-rank and discretization errors, which we do by developing a spatially accurate augmented BUG integrator. The novelties of the work can be summarized as follows,
\begin{itemize}
    \item \textit{A rank-adaptive low-rank Monte Carlo (RaDLR-MC) estimator}: We propose a rank-adaptive RaDLR-MC estimator based on the augmented BUG integrator, in which the rank of the solution is selected automatically at each time step to satisfy a prescribed truncation tolerance. Unlike the fixed-rank estimator of \cite{patwardhan_low-rank_2026}, this allows the approximation to adapt to evolving dynamics of the problem with reduced computational overhead.
    \item \textit{A novel rank-adaptive low-rank multilevel Monte Carlo (RaDLR-MLMC) framework}: We propose a low-rank MLMC framework based on geometric spatial refinement, a fixed discretization in velocity, and the augmented BUG integrator~\cite{ceruti_rank-adaptive_2022} for time-stepping on the low-rank manifold. The rank-adaptivity of the augmented BUG integrator is designed to preserve spatial discretization errors on each level, thereby ensuring that the MLMC theorem holds.
    \item \textit{Numerical experiment for high-dimensional kinetic equations}: We validate the proposed RaDLR-MC and RaDLR-MLMC estimators on benchmark kinetic equations, including radiative transport, radiation therapy, and shallow water moment equations. 
\end{itemize}

The rest of the paper is organized as follows: In section~\ref{sec:DLRA}, we recap DLRA for deterministic PDEs and under Monte Carlo sampling. We present a rank-adaptive low-rank Monte Carlo (RaDLR-MC) estimator for kinetic equations and discuss the relation between the discretization of the phase space and the rank of the approximation. In section~\ref{sec:lrMLMC}, we present a rank-adaptation strategy to obtain an augmented BUG integrator that has dominant spatial error. Based on this, we propose a new low-rank Multilevel Monte Carlo (RaDLR-MLMC) framework for kinetic equations. Finally, in section~\ref{sec:NumercalExperiments} we present numerical results for radiation transport, radiation therapy, and shallow water moment equations under uncertainty.

\section{Dynamical low-rank approximation and Monte Carlo sampling}\label{sec:DLRA}

\subsection{Fundamentals of DLRA for deterministic problems}
We begin this section by presenting the fundamental ideas of the continuous formulation of DLRA~\cite{koch_dynamical_2007} for deterministic time-dependent problems then extend it to the uncertain PDE case. 

Let $u(t,x,v)$ denote the solution to the deterministic version of \eqref{eq:genericKEwU}:
\begin{equation}\label{eq:genericKEwoU}
     \begin{aligned}
        \partial_{t}u(t,x,v) + \mathcal{A}(u)\partial_{x}u(t,x,v)  &= \mathcal{S}(u)(t,x,v),\\
        u(0,x,v) &= u_{0}(x,v),\\
        u(t,x,v) &= u_{b}(t,x,v),\quad x\in\Gamma_{\mathcal{D}}.
    \end{aligned}
\end{equation}
We restrict the solution of the deterministic kinetic equation to the low-rank manifold 
\begin{displaymath}
    \begin{aligned}
        \mathcal{M}_{r} :=\left\{g\in L^{2}(\mathcal{D}\times\mathcal{V}) \,|\,\right.& g(x,v) = \sum_{i,j=1}^{r}X_{i}(x)S_{ij}W_{j}(v)\\ &\left. \text{ where } X_{i}\in L^{2}(\mathcal{D}), W_{j}\in L^{2}(\mathcal{V}), S_{ij}\in \mathbb{R}\right\}
    \end{aligned}
\end{displaymath} 
such that $u$ admits the ansatz:
\begin{equation}\label{eq:lowrankAnsatz}
    u(t,x,v) \approx u_{r}(t,x,v) = \sum_{i,j=1}^{r}X_{i}(t,x)S_{ij}(t)W_{j}(t,v),
\end{equation}
where $\{X_{i}\}_{i=1,\ldots,r}\in L^{2}(\mathcal{D})$ and $\{W_{j}\}_{j=1,\ldots,r}\in L^{2}(\mathcal{V})$ are orthonormal bases, while $S_{ij}(t)$ denotes the mixing coefficients. Let $\langle \cdot,\cdot\rangle_{x}$ and $\langle \cdot,\cdot\rangle_{v}$ denote the inner products defined on $L^{2}(\mathcal{D})$ and $L^{2}(\mathcal{V})$, respectively. Then, by definition, we have
\begin{displaymath}
    \langle X_{i},X_{k}\rangle_{x} = \delta_{ik}, \quad \langle W_{j},W_{l}\rangle_{v} = \delta_{jl}\,.
\end{displaymath}
The goal is to evolve $X_{i}$, $W_{j}$, and $S_{ij}$ in time, such that the solution remains on the low-rank manifold $\mathcal{M}_{r}$ and follows the dynamics of the problem. In DLRA~\cite{koch_dynamical_2007}, this is done by projecting the dynamics of the problem onto the tangent space to the low-rank manifold. 
That is, if $\mathcal{P}_{u_{r}}$ denotes the orthogonal projection onto the tangent space to the low-rank manifold $\mathcal{M}_{r}$ at $u_{r}$, then we solve
\begin{equation}\label{eq:ProjectedDLRAeqn}
    \partial_{t}u_{r} = \mathcal{P}_{u_{r}}\left( \mathcal{F}(t,u_{r}) \right),
\end{equation}
where 
\begin{equation}\label{eq:RHSOperatorForm}
    \mathcal{F}(t,u_{r}) := -\mathcal{A}(u_{r})\partial_{x}u_{r} + \mathcal{S}(t,u_{r})\,.
\end{equation}
Using the short-hand notation $\mathcal{F}:=\mathcal{F}(t,u_{r})$ the orthogonal projection is given by
\begin{displaymath}
    \mathcal{P}_{u_{r}}\left( \mathcal{F}\right) = \sum_{j = 1}^{r}\left\langle W_{j},\mathcal{F} \right\rangle_{v}W_{j} - \sum_{i,j=1}^{r}X_{i}\left\langle X_{i}W_{j},\mathcal{F} \right\rangle_{x,v}W_{j} + \sum_{i = 1}^{r}X_{i}\left\langle X_{i},\mathcal{F} \right\rangle_{x}.
\end{displaymath}
By imposing the gauge conditions $\langle \partial_{t}X_{i},X_{j} \rangle_{x} = 0$ , $\langle \partial_{t}W_{i},W_{j} \rangle_{v} = 0$ we obtain unique evolution equations for $X_{i}$, $S_{ij}$, and $W_{j}$. Note that the specific gauge conditions have been chosen to obtain unique factorization of the low-rank factors in the tangent space~\cite{koch_dynamical_2007}. 

Since the rank of the solution is not known beforehand, it is often over-approximated, resulting in evolution equations for the factors that are then severely ill-conditioned. This is due to the inversion of the coefficient matrix $\mathbf{S}(t) = (S_{ij}(t))$, which has near-zero singular values, that is required for determining the evolution of the bases~\cite{koch_dynamical_2007}. Several robust integrators have been proposed in recent years that are robust to the presence of small singular values. The projector-splitting integrator (PSI), based on a Lie-Trotter or Strang splitting of the orthogonal projection, was introduced in \cite{lubich_projector-splitting_2014} and proven to be robust to the presence of small singular values~\cite{kieri_discretized_2016}. However, the PSI involves a backward in time step to update $S_{ij}(t)$, which leads to unstable dynamics for parabolic problems and even for hyperbolic problems~\cite{kusch_stability_2023}. The fixed-rank BUG integrator~\cite{ceruti_unconventional_2022}, also known as the unconventional integrator, was developed to mitigate this drawback. This development has facilitated the widespread use of DLRA integrators for kinetic problems~\cite{einkemmer_review_2025}. Since its development, the BUG integrator has been modified to allow rank-adaptivity~\cite{ceruti_rank-adaptive_2022}, completely parallel update steps~\cite{ceruti_parallel_2024}, and provably high-order integrators~\cite{ceruti_robust_2024,kusch_second-order_2025,nobile_high-order_2026}. In this work, we use the augmented BUG integrator~\cite{ceruti_rank-adaptive_2022}, a first-order in time, rank-adaptive integrator for time-dependent problems. 

We begin by deriving the evolution equations underlying all the BUG integrators. Let $K_{j}(t,x) := \sum_{i=1}^{r}X_{i}(t,x)S_{ij}(t)$ such that $u_{r}(t,x,v) = \sum_{j=1}^{r}K_{j}(t,x)W_{j}(t,v)$. Differentiating with respect to $t$ yields
\begin{displaymath}
    \partial_{t}u_{r} = \sum_{j=1}^{r}\partial_{t}\left( K_{j}W_{j} \right) = \sum_{j=1}^{r}\left( \partial_{t}K_{j}W_{j} + K_{j}\partial_{t}W_{j} \right).
\end{displaymath}
Inserting the above expression into the projected equation~\eqref{eq:ProjectedDLRAeqn} and projecting onto the space spanned by $\{W_{j}\}_{j=1,\ldots,r}$. For $1\leq l\leq r$, we obtain
\begin{displaymath}
    \begin{split}
        \sum_{j=1}^{r} \left(\partial_{t}K_{j}\langle W_{j}, W_{l}\rangle_{v} + K_{j}\langle \partial_{t}W_{j}, W_{l}\rangle_{v}\right) &= \sum_{j = 1}^{r}\left\langle W_{j},\mathcal{F} \right\rangle_{v}\langle W_{j},W_{l}\rangle_{v} \\&\hspace{-2cm}- \sum_{i,j=1}^{r}X_{i}\left\langle X_{i}W_{j},\mathcal{F} \right\rangle_{x,v}\langle W_{j},W_{l}\rangle + \sum_{i = 1}^{r}X_{i}\left\langle X_{i}W_{l},\mathcal{F} \right\rangle_{x,v}\,.
    \end{split}
\end{displaymath}
Using the orthogonality of $W_{j}$ and the gauge condition $\langle \partial_{t}W_{i},W_{j} \rangle_{v} = 0$, we arrive at
\begin{displaymath}
    \partial_{t}K_{l}(t,x) = \langle W_{l},\mathcal{F} \rangle_{v}.
\end{displaymath}
Note that since $W_{l}$ on the right-hand side is also time-dependent, we need to close the system. The system is closed by assuming that $W_{l}$ remains constant over the time step. Thus, to update the spatial basis from $t_{0}$ to $t_{1} = t_{0} + \Delta t$, we solve the following system of equations
\begin{equation}\label{eq:K_step_cont}
        \partial_{t}K_{l}(t,x) = \langle W_{l},\mathcal{F} \rangle_{v}\,, \quad \partial_{t}W_{l} = 0.
    \end{equation}

Similarly, defining $L_{i}(t,v) := \sum_{j=1}^{r}S_{ij}(t)W_{j}(t,v)$ and following the same steps as for $K_{j}$, we obtain the following evolution equations for the remaining factors:
\begin{align}
    \partial_{t}L_{i}(t,v) &= \langle X_{i},\mathcal{F} \rangle_{x},\quad \partial_{t}X_{i}(t,x) = 0,\label{eq:L_step_cont}\\
    \partial_{t}S_{ij}(t) &= \langle X_{i}W_{j},\mathcal{F} \rangle_{x,v}\quad \partial_{t}X_{i}(t,x) = 0,\,\partial_{t}W_{j}(t,v) = 0\,. \label{eq:S_step_cont}
\end{align}
The update for $S_{ij}$ is obtained by testing \eqref{eq:ProjectedDLRAeqn} with $X_{i}W_{j}$ and integrating over the spatial and velocity domains. The choice of initial conditions for these closure equations leads to various integrators. A summary of the integrators is provided in Table \ref{tab:Comparison_low_rank_integrators}. The initial conditions for the K- and L-steps of all the integrators are the same as for the augmented BUG integrator described below. 
\begin{table}
    \renewcommand{\arraystretch}{1.4}
    \small
    \begin{tabular}{lccc}
        \toprule
        Integrator & $S(t_0)$ & rank-adaptive & parallel updates \\
        \midrule
        fixed-rank BUG~\cite{ceruti_unconventional_2022}
            & $X^{1,\top} X^0 S^0 W^{0,\top} W^{1}\in\R^{r\times r}$
            & no & partial \\
        fixed-rank augmented BUG~\cite{ceruti_rank-adaptive_2022}
            & $\widetilde{X}^\top X^{0} S^{0} W^{0,\top}\widetilde{W}\in\R^{2r\times 2r}$
            & no & partial \\
        augmented BUG~\cite{ceruti_rank-adaptive_2022}
            & $\widetilde{X}^\top X^{0} S^{0} W^{0,\top} \widetilde{W}\in\R^{2r\times 2r}$
            & yes & partial \\
        parallel BUG~\cite{ceruti_parallel_2024}
            & $S_0\in\R^{r\times r}$
            & yes & full \\
        \bottomrule
    \end{tabular}
    \caption{Comparison of BUG integrators for dynamical low-rank approximation, including the S-step initial condition, rank-adaptivity, and parallelism properties. $\widetilde{X}$ and $\widetilde{W}$ are as defined in the K- and L-steps of the augmented BUG integrator. ``Partial" parallel updates are in the K- and L-steps, while ``full" implies in all three sub-steps.}
    \label{tab:Comparison_low_rank_integrators}
\end{table}

In the augmented BUG integrator, \eqref{eq:K_step_cont} and \eqref{eq:L_step_cont} are first updated in parallel from initial data. The basis is expanded by augmenting the previous basis at time $t_{0}$ to the updated basis. The coefficients are subsequently updated by projecting onto the augmented basis and evolving time with \eqref{eq:S_step_cont}. To make this more precise, we present one step of the augmented BUG integrator. For the given initial data $u_{r}(t_{0},x,v) = \sum_{i,j=1}^{r_{0}}X_{i}^{0}(x)S_{ij}^{0}W_{j}^{0}(v) = X^{0,\top}(x)\mathbf{S}^{0}V^{0}(v)$ at time $t_{0}$ the augmented BUG integrator computes the solution at time $t_{1} = t_{0} + \Delta t$ as follows:

\begin{enumerate}
    \item \textbf{K-step}: Update and augment the spatial basis from $\{X^{0}_{i}(x)\}_{i=1}^{r_{0}}$ to $\{\widetilde{X}^{1}_{i}\}_{i=1}^{2r_{0}}$ by solving
    \begin{displaymath}
        \partial_{t}K_{j} = \langle W_{j}^{0}, \mathcal{F}_{K} \rangle_{v}, \qquad K_{j}(t_{0},x) = \sum_{i=1}^{r_{0}}X_{i}^{0}(x)S_{ij}^{0},
    \end{displaymath}
    where $\mathcal{F}_{K}(t,x,v) := \mathcal{F}(t,\sum_{j=1}^{r}K_{j}(t,x)W_{j}^{0}(v))$. Compute $\{ \widetilde{X}_{i}^{1}(x) \}_{i = 1,\ldots,2r_{0}} $ as an orthonormal basis of $\left\{ K_{1}^{1}(x),\ldots,K_{r_{0}}^{1}(x),X_{1}^{0}(x),\ldots,X_{r_{0}}^{0}(x)\right\}$ and store the matrix $\mathbf{M}\in\mathbb{R}^{2r_{0}\times r_{0}}$ with entries $(\mathbf{M})_{\alpha\beta} = \langle \widetilde{X}^{1}_{\alpha},X^{0}_{\beta}\rangle_{x}$.
    \item \textbf{L-step}: Update and augment the velocity basis from $\{W^{0}_{j}(v)\}_{i=1}^{r_{0}}$ to $\{\widetilde{W}^{1}_{j}\}_{j=1}^{2r_{0}}$ by solving
    \begin{displaymath}
        \partial_{t}L_{i} = \langle X_{i}^{0},\mathcal{F}_{L} \rangle_{x}, \qquad L_{i}^{0}(v) = \sum_{j=1}^{r_{0}}S_{ij}^{0}W_{j}^{0}(v),
    \end{displaymath}
    where $\mathcal{F}_{L}(t,x,v) :=\mathcal{F}(t,\sum_{i=1}^{r}X_{i}^{0}(x)L_{i}(t,v)) $. Compute $\{ \widetilde{W}_{j}^{1} \}_{j = 1,\ldots,2r_{0}} $ as an orthonormal basis of $\left\{ L_{1}^{1},\ldots,L_{r_{0}}^{1},W_{1}^{0},\ldots,W_{r_{0}}^{0} \right\}$ and store the matrix $\mathbf{N}\in\mathbb{R}^{2r_{0}\times r_{0}}$ with entries $(\mathbf{N})_{\alpha\beta} = \langle \widetilde{W}^{1}_{\alpha},W^{0}_{\beta}\rangle_{v}$.
    \item Define $\widetilde{\mathbf{S}}^{0} = \mathbf{M}\mathbf{S}^{0}\mathbf{N}^{\top}\in\mathbb{R}^{2r_{0}\times 2r_{0}}$, where $\mathbf{S}^{0} = (S_{ij}^{0})\in\mathbb{R}^{r_{0}\times r_{0}}$. Update the coefficients to $\widetilde{\mathbf{S}}^{1} = \widetilde{\mathbf{S}}(t_{1})\in\mathbb{R}^{2r_{0}\times2r_{0}}$ by solving
    \begin{displaymath}
        \frac{d}{dt}\widetilde{S}_{ij}(t) = \langle \widetilde{X}_{i}^{1}\widetilde{W}_{j}^{1}, \mathcal{F}_{S}\rangle_{x,v}, \qquad \widetilde{S}_{ij}(t_{0}) = (\widetilde{\mathbf{S}}^{0})_{ij},
    \end{displaymath}
    where $\mathcal{F}_{S}(t,x,v) := \mathcal{F}(t,\sum_{i,j=1}^{r}\widetilde{X}_{j}^{1}(x)\widetilde{S}_{ij}(t)\widetilde{W}_{j}^{1}(v))$.
    \item \textbf{Truncation}: Compute the SVD $\widetilde{\mathbf{S}} = \widetilde{\mathbf{P}}\widetilde{\boldsymbol{\Sigma}}\widetilde{\mathbf{T}}^{\top}$, where $\boldsymbol{\Sigma} = \mathrm{diag}(\sigma_{k})$ contains the singular values of $\widetilde{\mathbf{S}}$. Choose the new rank $r_{1}\leq2r_{0}$ as the smallest number satisfying:
        \begin{displaymath}
            \left( \sum_{k = r_{1}+1}^{2r_{0}}\sigma_{k}^{2} \right)^{1/2}\leq \vartheta\,.
        \end{displaymath}
         For $i,j = 1,\ldots,r_{1}$, set $S_{ij}^{1} = \sigma_{i}\delta_{ij}$, $X_{i}^{1}(x) = \sum_{k=1}^{2r_{0}}\widetilde{X}_{k}^{1}(x)\widetilde{P}_{ki}$, and $W_{j}^{1}(v) = \sum_{k=1}^{2r_{0}}\widetilde{W}_{k}^{1}(v)\widetilde{T}_{kj}$ 
\end{enumerate}
Then the factorized solution at time $t_{1}$ with rank $r_{1}$ is given by $\sum_{i,j=1}^{r_{1}}X_{i}^{1}(x)S_{ij}^{1}V_{j}^{1}(v)$.

\begin{remark}
     Note that in practice, to update the low-rank factors, we need to define an appropriate discretization of the spatial and velocity domains and the corresponding inner products.
\end{remark}

\subsection{Rank-adaptive dynamical low-rank Monte Carlo (RaDLR-MC) estimator}

In the probabilistic setting, we treat the random variable as a parameter of the kinetic equation and assume that the solution to \eqref{eq:genericKEwU} admits the following uncertainty-dependent low-rank ansatz:
\begin{equation}\label{eq:randomLowRankAnsatz}
    u(t,x,v\,;\omega) \approx u_{r}(t,x,v\,;\omega) = \sum_{i,j=1}^{r}X_{i}(t,x\,;\omega)S_{ij}(t\,;\omega)W_{j}(t,v\,;\omega)\,.
\end{equation}
 Thus, for each realization of the random variable, we can compute the low-rank approximation by repeatedly applying the augmented BUG integrator with a given truncation tolerance $\vartheta$. This estimator can also be extended to any (rank-adaptive) BUG integrator with a pre-defined truncation tolerance $\vartheta$ in a straightforward way. Let $Q_{r}$ denote the low-rank approximation to the QoI $Q$, computed with the augmented BUG integrator. To estimate $\mathbb{E}\left[Q\right]$ defined in \eqref{eq:ExpectedVal}, we construct and compute the estimator $\widehat{Q}_{r}$ to $\mathbb{E}\left[Q_{r}\right]$. For completeness, we state the RaDLR-MC estimator. Let $\omega_{1},\ldots,\omega_{M}$ be $M$ independent and identically distributed (i.i.d.) random realizations of $\omega$. Then, we define the RaDLR-MC estimator as
\begin{displaymath}
    \widehat{Q}_{r}^{MC} := \frac{1}{M}\sum_{i=1}^{M}Q_{r}^{(i)}
\end{displaymath}
where $Q_{r}^{(i)} := Q_{r}(\omega_{i})$ denotes the QoI approximated with the augmented BUG integrator with truncation tolerance $\vartheta$ for the $i^{th}$ i.i.d. copy of $\omega$. We quantify the accuracy of the estimator with the mean squared error (MSE) defined as
\begin{displaymath}
    e\left(\widehat{Q}_{r}^{MC}\right)^{2} \coloneqq \mathbb{E}\left[\left\lVert \widehat{Q}_{r} - \mathbb{E}\left[Q\right] \right\rVert_{H}^{2} \right]. 
\end{displaymath}
The estimator's dependence on the truncation tolerance $\vartheta$ can be made more explicit by replacing the subscript $r$ by $\vartheta$ in $u_{r}$ and $Q_{r}$. However, to avoid introducing new notation, we continue with the subscript $r$ to denote the low-rank approximation. In the context of rank-adaptive integrators, $r$ should be interpreted as a varying quantity that depends on $\vartheta$. 

 Unlike the fixed-rank augmented BUG integrator in \cite{patwardhan_low-rank_2026}, the rank in the augmented BUG integrator evolves as a function of the truncation tolerance $\vartheta$. Although this approach improves adaptability to the underlying problem dynamics and uncertainty, selecting an appropriate value for $\vartheta$ demands a thorough understanding of the integrator, the problem at hand, and the numerical method employed. For instance, the parallel BUG integrator requires a smaller tolerance compared to the augmented BUG integrator to get comparable results for the same test case~\cite{ceruti_parallel_2024}. Moreover, if $R(t\,;u_{r},\vartheta)$ denotes the rank of the approximation $u_{r}$ at time $t$ for truncation tolerance $\vartheta$, then in the probabilistic setting $R$ is a random variable. Thus, the error in the low-rank approximation may vary from highly accurate to highly inaccurate across samples. 

\subsubsection{Evaluating the bias of augmented BUG}\label{sec:DiscLRA}
To gain more insight into the relation between the low-rank approximation and truncation tolerance, we evaluate the bias of the RaDLR-MC estimator. The MSE of the RaDLR-MC estimator can be decomposed as 
\begin{displaymath}
    e\left(\widehat{Q}_{r}^{MC}\right)^{2} = \mathbb{E}\left[\left\lVert \frac{1}{M}\sum_{i=1}^{M}(Q_{r}^{(i)} - \mathbb{E}[Q_{r}]) \right\rVert^{2}_{H}\right] + \lVert\mathbb{E}\left[Q_{r} - Q\right] \rVert_{H}^{2}\,,
\end{displaymath}
where we have used the fact that $H$ is a Hilbert space with the inner product norm $\lVert \cdot\rVert^{2}_{H} = \langle\cdot,\cdot\rangle_{H}$ and 
\begin{displaymath}
    \begin{aligned}
        \mathbb{E}\left[\left\langle\widehat{Q}_{r}^{MC} - \mathbb{E}\left[Q_{r}\right],\mathbb{E}\left[Q_{r}\right] - \mathbb{E}\left[Q\right]  \right\rangle_{H}\right] &= \langle\, \underset{= \mathbb{E}[Q_{r}]}{\underbrace{\mathbb{E}[\widehat{Q}_{r}^{MC}]}} - \mathbb{E}\left[Q_{r}\right],\mathbb{E}\left[Q_{r}\right] - \mathbb{E}\left[Q\right] \,\rangle_{H} = 0
    \end{aligned}
\end{displaymath}
Thus, the MSE can be written in terms of the estimator variance $V_{r}$ and bias $B_{r}$ as
\begin{displaymath}
    e\left(\widehat{Q}_{r}^{MC}\right)^{2} = \frac{1}{M}V_{r} + B_{r}^{2},
\end{displaymath}
where 
\begin{displaymath}
    V_{r} := \mathbb{E}[\lVert Q_{r}^{(1)} - \mathbb{E}[Q_{r}] \rVert^{2}_{H}] \quad \mathrm{and}\quad B_{r} := \lVert \mathbb{E}[Q_{r} - Q] \rVert_{H}\,.
\end{displaymath}

We can further decompose the bias $B_{r}$ into contributions from the spatio-velocity discretization and the low-rank approximation. Let $u_{full}$ denote the approximation to \eqref{eq:genericKEwU} computed using the same numerical discretization method (without DLRA) as for the low-rank approximation $u_{r}$, and let $Q_{full}$ denote the corresponding QoI. We assume that $u_{full}$ and $u_{r}$ have the exact same spatio-velocity discretization parameters, i.e., $\Delta x_{full} = \Delta x = \Delta x_{r}$ and $\Delta v_{full} = \Delta v = \Delta v_{r}$, where $\Delta x\in\mathbb{R}^{d}_{>0}$ and $\Delta v\in\mathbb{R}^{d}_{>0}$ are the spatial and velocity discretization parameters. Adding and subtracting $Q_{full}$ in the bias and using the triangle inequality, we obtain
\begin{displaymath}
    B_{r} \leq \underset{\text{low-rank contr.}}{\underbrace{\lVert\mathbb{E}\left[Q_{r} - Q_{full}\right] \rVert_{H}}} + \underset{\text{spatio-velocity contr.}}{\underbrace{\lVert\mathbb{E}\left[Q_{full} - Q\right] \rVert_{H}}}\,.
\end{displaymath}
To further simplify the inequality, we use Jensen's inequality and the Lipschitz continuity of $Q$, which yields
\begin{displaymath}
    B_{r} \leq L\cdot\mathbb{E}[\lVert u_{r} - u_{full} \rVert] + L\cdot\mathbb{E}[\lVert u_{full} - u \rVert]\,,
\end{displaymath}
where $\lVert\cdot\rVert$ is an appropriate norm for the approximations $u_{r}$ and $u_{full}$ to \eqref{eq:genericKEwU}, for instance the $L^{2}$ norm.

The error contribution of the spatio-velocity discretization can be estimated using standard theory on numerical methods. For an order $\alpha\geq 1$ method in space and order $q\geq 1$ in velocity domain, we assume that the error contribution of the spatio-velocity discretization has the form 
\begin{displaymath}
    \lVert u - u_{full} \rVert \leq K_{1}(\omega)h_{x}^{\alpha} + K_{2}(\omega)h_{v}^{q}\,,
\end{displaymath}
where $h_{x} := \lVert \Delta x\rVert_{\infty}$, $h_{v} := \lVert \Delta v\rVert_{\infty}$, and $\lVert \cdot\rVert_{\infty}$ is the infinity norm.

In \cite{ceruti_rank-adaptive_2022}, the authors prove rigorous error bounds for the augmented BUG integrator that are robust to the presence of small singular values (over-approximation of the rank). However, these error bounds use the Lipschitz continuity of the right-hand side operator $\mathcal{F}$. For kinetic equations, the Lipschitz constant of the right-hand side scales as $\mathcal{O}(1/h_{x})$, rendering the error bounds invalid.

While no rigorous bounds exist for the augmented BUG integrator without the Lipschitz constant, from experiments on deterministic problems, it is known that the error tends to respect the following bound
\begin{displaymath}
    \lVert u_{r}(t_{s},\cdot,\cdot) - u_{full}(t_{s},\cdot,\cdot)\rVert \leq K_{3}\Delta t + K_{4}s\vartheta + K_{5}\epsilon,
\end{displaymath}
where $t_{s} = t_{0} + s\Delta t$ and $\Delta t = t_{\mathrm{end}}/N_{t}$ for some $N_{t}\in\mathbb{N}$, $\epsilon$ is the low-rank projection error, and $K_{3},\,K_{4},\,K_{5}$ are independent of the singular values of $u_{r}$. Overall, the bias term has the following structure
\begin{displaymath}
    B_{r} \leq L(\widehat{K}_{1}h_{x}^{\alpha} + \widehat{K}_{2}h_{v}^{q} + \widehat{K}_{3}\Delta t + \widehat{K}_{4}s\vartheta+ \widehat{K}_{5}\epsilon),
\end{displaymath}
where $L$ is the Lipschitz constant of $Q$ and $\widehat{K}_{i}$, $i\in\{1,2,3,4,5\}$, are expectations of $K_{i}$ over $p(\omega)$. 

\subsubsection{Cost analysis of the RaDLR-MC estimator}
\label{sec:costAnalysis}
We evaluate the cost of achieving an MSE of $\varepsilon^{2}$ for the RaDLR-MC estimator, i.e.,
\begin{align*}
    e\left(\widehat{Q}_{r}^{MC}\right)^{2} &\leq \varepsilon^{2},\\
    \implies \frac{1}{M}V_{r} + B_{r}^{2} &\leq \varepsilon^{2}.
\end{align*}
For $0<\theta< 1 $, the contribution of the statistical error and numerical error can be distributed as 
\begin{displaymath}
    \frac{1}{M}V_{r} \leq \theta\cdot\varepsilon^{2}, \qquad B_{r} \leq \sqrt{1-\theta}\cdot\varepsilon\,.
\end{displaymath}
This splitting, of the contributions to the MSE, leads to a trade-off between the accuracy of the approximation and the variance of the estimator. Usually, one is interested in computing a low-variance estimate of an accurate approximation of the QoI. In terms of the low-rank approximation error, this implies that the solution is computed on a fine grid in the phase space with a small truncation tolerance to accurately capture the solution. Moreover, $M\geq V_{r}/(\theta\cdot\varepsilon^{2})$ samples must be performed to achieve the required low variance. Thus to construct such an estimator, we need to completely understand the relation between the bias and its various contributions. The behavior of the approximation error due to spatial and velocity discretizations is a well studied subject and can be found in any standard textbook on the topic. See \cite{leveque_finite_2002,case_linear_1967} for instance. A less understood relation, however, is that between the approximation error and the truncation tolerance $\vartheta$. To study this for the augmented BUG integrator, we consider the linesource test case, described in \cite{ceruti_rank-adaptive_2022}.

\begin{figure}
    \centering
    \begin{subfigure}[b]{0.45\linewidth}
        \includegraphics[width=\linewidth]{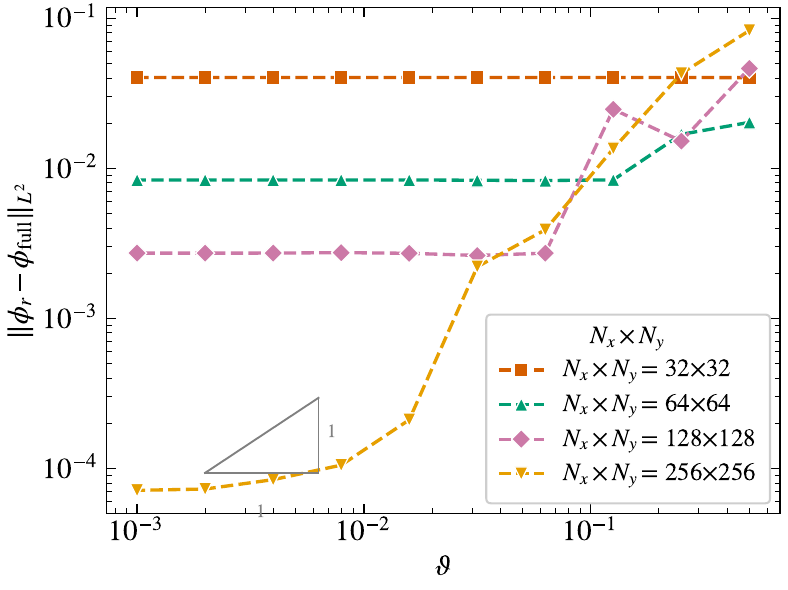}
    \end{subfigure}
    \begin{subfigure}[b]{0.45\linewidth}
        \includegraphics[width=\linewidth]{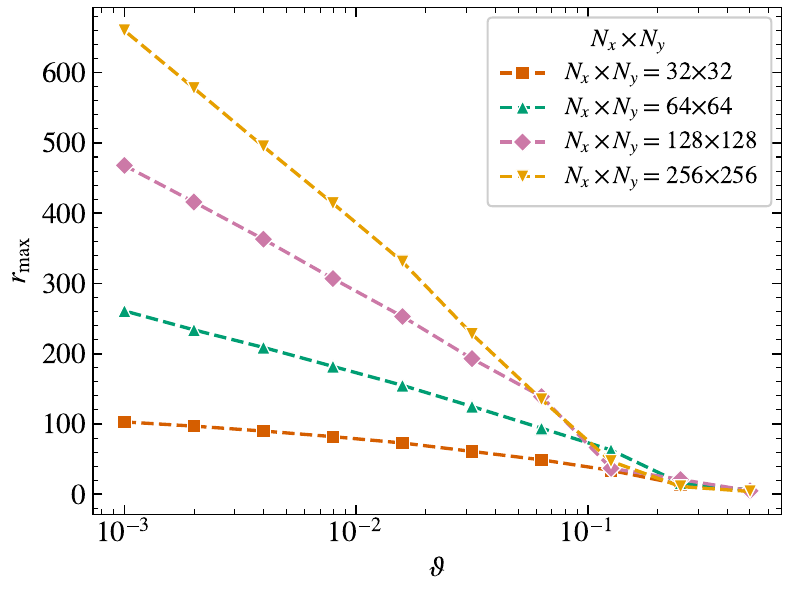}
    \end{subfigure}

    \begin{subfigure}[b]{0.45\linewidth}
        \includegraphics[width=\linewidth]{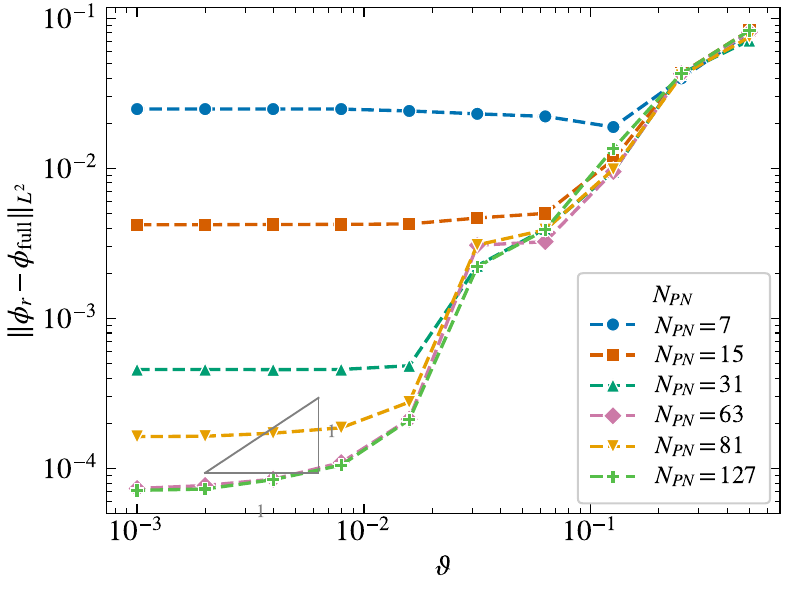}
    \end{subfigure}
    \begin{subfigure}[b]{0.45\linewidth}
        \includegraphics[width=\linewidth]{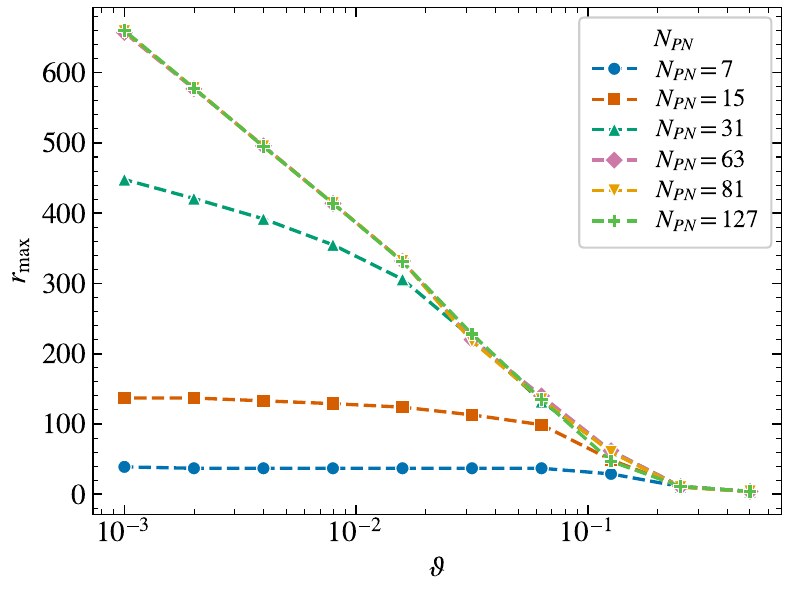}
    \end{subfigure}
    \caption{Dependence of the low-rank approximation error and maximum rank on the truncation tolerance $\vartheta$ for the line source test case~\cite{ceruti_rank-adaptive_2022}. Here $Q(x) := \phi(t=1.0,x) = \int_{S^2} \psi(t=1.0,x,\Omega)\, \mathrm{d}\Omega$ denotes the scalar flux and $\Omega = v/\lvert v\rvert$ is the direction of flight. (Top) Spatial refinement study: discrete $L^2$ error of $\phi_r$ against a full-rank reference $\phi_\mathrm{full}$ on a $256\times256$ grid (left) and maximum rank $r_{\max}$ over all time steps (right), for spatial grids ranging from $32\times32$ to $256\times256$ cells with fixed angular resolution. (Bottom) Angular refinement study: same quantities for angular resolutions ranging from $N=8$ to $N=128$ with fixed spatial grid.}
    \label{fig:Error_vs_tolerance}
\end{figure}

In Figure \ref{fig:Error_vs_tolerance}, we plot the absolute error of the low-rank approximation of the QoI, $Q_{r}$, to the full rank approximation $Q_{full}$ against the truncation tolerance of the augmented BUG integrator for various spatial grids and velocity discretizations. 

From Figure~\ref{fig:Error_vs_tolerance} we see that, the error depends on both the spatio-velocity discretization and the truncation tolerance. Moreover, for a fixed spatial and velocity discretization, it is not possible to achieve a target MSE by varying the truncation tolerance alone (and hence the approximation rank), unless the spatio-velocity error is already bounded by $\varepsilon/\sqrt{2}$. This limitation extends to the low-rank control variates strategy proposed in \cite{patwardhan_low-rank_2026}. Nevertheless, Figure~\ref{fig:Error_vs_tolerance} shows that refining the spatial grid at a fixed velocity resolution reduces the approximation error, and thus the bias and overall MSE of the estimator. We exploit this observation in the next section to construct a spatially accurate augmented BUG integrator and an adaptive estimator for kinetic equations.

\section{A rank-adaptive framework for low-rank multilevel Monte Carlo}\label{sec:lrMLMC}
Exploiting spatial refinement as the primary mechanism for controlling approximation error, we construct a multilevel estimator over a hierarchy of spatial discretizations, each paired with a low-rank approximation of the kinetic equation. The intuition is straightforward: coarse levels are cheap to evaluate and capture the bulk of the solution structure at low rank, while successive refinements incrementally correct the bias introduced by both spatial discretization and rank truncation. The multilevel telescoping sum then aggregates these corrections, and the sample allocation across levels can be tuned to balance statistical and approximation errors. The outcome is an estimator that pairs the well-known variance reduction of MLMC with the computational savings of low-rank compression in phase space; a combination particularly well-suited to the high-dimensional setting of kinetic equations. We call this the \textit{rank-adaptive dynamical low-rank multilevel Monte Carlo} (RaDLR-MLMC) estimator.

\subsection{A spatially accurate low-rank approximation}
The error of the low-rank approximation of the QoI is governed by four discretization parameters: spatial cell size $\Delta x$, velocity grid cell size $\Delta v$, time step size $\Delta t$, and truncation tolerance $\vartheta$ (recall Figure~\ref{fig:Error_vs_tolerance}). To lay the groundwork for the multilevel estimator, we first describe how spatial refinement can be incorporated into the low-rank approximation framework.

We use an explicit Euler time-stepping scheme and for stability considerations, the time step size is chosen to satisfy the following Courant-Friedrichs-Lewy (CFL) condition
\begin{equation}
    \Delta t \leq \frac{\mathrm{CFL}}{\lvert v\rvert\sigma}\left(\dfrac{\prod_{i=1}^{d}\Delta x_{i}}{\sum_{i=1}^{d}\prod_{i\neq j}\Delta x_{j}}\right),
\end{equation}
where $\lvert v\rvert := \sup(\mathcal{V})<\infty$. Note that the CFL condition may contain additional terms, like the scattering or absorption coefficient for radiation transport. These problem-specific terms are incorporated in $\sigma$. On a uniform spatial grid $\Delta x_{i} = h_{x}$, for $i = 1,\ldots,d$, the product $\prod_{i=1}^{d}\Delta x_{i} = h_{x}^{d}$ and the sum of products is $\sum_{i=1}^{d}\prod_{i\neq j}\Delta x_{j} = d\,h^{d-1}_{x}$, so the CFL condition reduces to
\begin{equation}
    \Delta t \leq \dfrac{\mathrm{CFL}}{d\,\lvert v\rvert\sigma}h_{x}\,.
\end{equation}
For an order $\alpha$ method in space, choosing the truncation tolerance to satisfy
\begin{equation}
    \vartheta\leq \frac{\widehat{K}_{1}}{\widehat{K}_{4}}\left( \dfrac{\mathrm{CFL}}{d\,\sigma\lvert v\rvert\, t_{end}} \right)h_{x}^{\bar{\alpha}},
\end{equation}
where $\bar{\alpha} := \min\{\alpha+1,2\}$, ensures that the low-rank approximation error does not dominate the spatial discretization error. Since we use an asymptotically first order method in time, using a higher-order spatial discretization does not yield additional benefits. This could be mitigated by using a higher-order BUG integrator~\cite{ceruti_robust_2024,kusch_second-order_2025,nobile_high-order_2026} in conjunction with a higher-order spatial discretization. 

Substituting the CFL-constrained $\Delta t$ and the bound on $\vartheta$ into the general error estimate, and assuming a uniform spatial grid, the dominant term in the bias reduces to
\begin{equation}\label{eq:Bias_LR}
    B_{r} \leq \widehat{K}_{6} h_{x} + \widehat{K}_{2}h_{v}^{q} + \widehat{K}_{5}\epsilon\,,
\end{equation}
where
\begin{displaymath}
    \widehat{K}_{6} = L\cdot\max\left\{ \widehat{K}_{1}, \widehat{K}_{3},\widehat{K}_{1}\left(\dfrac{\mathrm{CFL}}{d\,\sigma\lvert v\rvert\, t_{end} }\right) \right\}\,.
\end{displaymath}

\subsection{The RaDLR-MLMC method and algorithm}
\label{sec:method}
Finally, we state the rank-adaptive dynamical low-rank multilevel Monte Carlo (RaDLR-MLMC) estimator for kinetic equations. We see from \eqref{eq:Bias_LR} that for a high-order method in velocity ($q>1$) and small low-rank projection error $\epsilon\ll1$, the bias is dominated by the spatial error. The assumption that the low-rank projection error is non-dominant remains valid for the problems considered in this work which require the solution to admit a low-rank structure. Thus, we use a high-order method in the velocity domain such that the spatial error remains the dominant contribution to the bias term and construct the multilevel hierarchy in space. Let $h_{x,\ell} \coloneqq h_{x,0}/2^{\ell}$ denote the width of the uniform spatial grids, where $\ell = 0,1,\ldots,L$ are referred to as levels and $h_{x,0}$ is the size of a spatial cell on the coarsest level. Furthermore, we let 
\begin{displaymath}
    \vartheta_{\ell} := \frac{\widehat{K}_{1}}{\widehat{K}_{4}}\left( \dfrac{\mathrm{CFL}}{d\,\sigma\lvert v\rvert\,t_{end}} \right)h_{x,\ell}^{\bar{\alpha}},\qquad \mathrm{and} \qquad \Delta t_{\ell} := \left(\dfrac{\mathrm{CFL}}{d\,\sigma\lvert v\rvert}\right)h_{x,\ell}\,.
\end{displaymath}
If we define the low-rank approximation on the grid parametrized by $(\Delta t_{\ell},h_{x,\ell},\Delta v,\vartheta_{\ell})$ by $Q_{\ell}$, $\ell=0,1,\ldots,L$, then $\mathbb{E}\left[Q_{L}\right]$, on the finest level $L$, can be written as the telescoping sum
\begin{equation}\label{eq:TelescopingSum}
    \mathbb{E}\left[Q_{L}\right] = \mathbb{E}\left[Q_{0}\right] + \sum_{\ell = 1}^{L}\left( \mathbb{E}\left[Q_{\ell}\right] - \mathbb{E}\left[Q_{\ell -1}\right] \right) = \sum_{\ell =0}^{L}\mathbb{E}\left[\Delta Q_{\ell}\right],
\end{equation}
where we define the random variable $\Delta Q_{\ell}(\omega) \coloneqq Q_{\ell}(\omega) - Q_{\ell-1}(\omega)$, with $Q_{-1}(\omega) \equiv 0$. The basic principle of the MLMC method is to estimate the QoI on the finest level $L$ by sampling from $\Delta Q_{\ell}$ on level $\ell = 0,\ldots,L$.

The RaDLR-MLMC estimator for $\mathbb{E}\left[Q_{L}\right]$ based on the telescoping sum~\eqref{eq:TelescopingSum} is defined as 
\begin{equation}\label{eq:MLMCEstimatorSTD}
    \widehat{Q}_{L}^{ML} \coloneqq \sum_{\ell=0}^{L}\frac{1}{M_{\ell}}\sum_{i = 1}^{M_{\ell}}\Delta Q_{\ell}^{(i)},
\end{equation}
where $M_{\ell}$ is the sample size on level $\ell$ and $\Delta Q_{\ell}^{(i)}$, $i = 1,\ldots,M_{\ell}$, are $M_{\ell}$ i.i.d. copies of $\Delta Q_{\ell}$. The MSE of the RaDLR-MLMC estimator is given by
\begin{displaymath}
         e\left(\widehat{Q}_{L}^{ML}\right)^{2} = \mathbb{E}\left[\left\lVert \widehat{Q}_{L}^{ML} - \mathbb{E}\left[Q\right] \right\rVert_{H}^{2} \right]= \mathbb{E}\left[\left\lVert \widehat{Q}_{L}^{ML} - \mathbb{E}\left[Q_{L}\right] \right\rVert_{H}^{2} \right] + \left\lVert \mathbb{E}\left[Q_{L} - Q\right] \right\rVert_{H}^{2}\,.
\end{displaymath}
Inserting the MLMC estimator~\eqref{eq:MLMCEstimatorSTD} and the telescoping sum~\eqref{eq:TelescopingSum} in the MSE, we obtain 
\begin{align*}
    e\left(\widehat{Q}_{L}^{ML}\right)^{2} &=  \mathbb{E}\left[\left\lVert \sum_{\ell=0}^{L}\frac{1}{M_{\ell}}\sum_{i = 1}^{M_{\ell}}\left(\Delta Q_{\ell}^{(i)} - \mathbb{E}\left[\Delta Q_{\ell}\right] \right) \right\rVert_{H}^{2} \right] + \left\lVert \mathbb{E}\left[Q_{L} - Q\right] \right\rVert_{H}^{2}\\
    &= \sum_{\ell=0}^{L}\frac{1}{M_{\ell}} \mathbb{E}\left[\left\lVert\Delta Q_{\ell} - \mathbb{E}\left[\Delta Q_{\ell}\right] \right\rVert_{H}^{2} \right] + \left\lVert \mathbb{E}\left[Q_{L} - Q\right] \right\rVert_{H}^{2}\\
    &= \sum_{\ell=0}^{L}\frac{1}{M_{\ell}}V_{\ell} + \left\lVert \mathbb{E}\left[Q_{L} - Q\right] \right\rVert_{H}^{2}
\end{align*}
where $V_{\ell} := \mathbb{E}\left[\left\lVert\Delta Q_{\ell} - \mathbb{E}\left[\Delta Q_{\ell}\right] \right\rVert_{H}^{2} \right] = \mathbb{V}[\Delta Q_{\ell}]$ has been defined in \eqref{eq:VarianceVal}. 

The first term of the MSE is called the total variance of the estimator, while the second term is the bias of the estimator. In the MLMC method, the statistical and numerical errors of the estimator are balanced to reach a certain predefined accuracy $\varepsilon$. To be specific, for a given $\varepsilon>0$, we want to choose the finest discretisation level $L$ and the number of samples $M_{\ell}$ on each level such that the MSE of the RaDLR-MLMC estimator satisfies $e\left(\widehat{Q}_{L}^{ML}\right)^{2} \leq \varepsilon^{2}$. As discussed in Section \ref{sec:costAnalysis}, this leads to a trade-off between the error contributions which can be distributed using a parameter $0<\theta\leq1$. One way to resolve this is by choosing $\theta=0.5$, i.e., equally splitting the $\varepsilon$-error between the two error contributions and thus requiring that the total variance and bias be at most $\varepsilon^{2}/2$. This splitting of the error is not necessarily the most optimal strategy~\cite{haji-ali_optimization_2016}, but rather a pragmatic compromise.

If $C_{\ell}$ is the cost of computing a single sample of $\Delta Q_{\ell}^{(i)}$, then the total cost of the RaDLR-MLMC estimator is given by
\begin{displaymath}
    \mathcal{C}\left( \widehat{Q}_{L}^{ML} \right) = \sum_{\ell=0}^{L}M_{\ell}C_{\ell}\,.
\end{displaymath}

If we have an a priori estimate of the bias, then we can find an optimal $L$ satisfying $\left\lVert \mathbb{E}\left[Q_{L} - Q\right] \right\rVert_{H}^{2}\leq \varepsilon^{2}/2$. In general, the constants are not known a priori, and hence $L$ is chosen on-the-fly by estimating the bias during computation as described in \cite{cliffe_multilevel_2011}.  

\begin{algorithm}
    \caption{RaDLR-MLMC algorithm}\label{alg:MLMC}
    \begin{algorithmic}[1]
        \State Initialize $L =2$ with initial samples $M_{\ell}$ on levels $\ell = 0,1,2$ 
        \State Initialize $K_{\ell} =0$ to be the count of already computed samples on level $\ell$
        \While{$K_{\ell} < M_{\ell}$}
        \For{$\ell=1,\ldots,L$}
          \State On level $\ell$, generate additional samples of $\Delta Q_{\ell}^{(i)}$ for $i = 1,\ldots,M_{\ell}-K_{\ell}$ with parameters $(\Delta t_{\ell}, h_{\ell}, \Delta v, \vartheta_{\ell})$ and the augmented BUG integrator
          \State Updates estimates of $V_{\ell}$ for $\ell = 0,\ldots,L$
          \State Compute the optimal number of samples $M_{\ell}$ on each level $\ell = 0,1,\ldots,L$ according to \eqref{eq:OptSamples}
          \If{Test for weak convergence of the QoI}
            \State If weak convergence fails, set $L\coloneqq L+1$, initialize $M_{L}$, and $K_{L}=0$
          \EndIf
          \EndFor
        \EndWhile
    \end{algorithmic}
\end{algorithm}

Assuming that we know the final level $L$, the number of samples on each level is chosen such that the total cost is minimized while the total variance is bounded by $\varepsilon^{2}/2$. That is, the optimal number of samples of level $\ell$ is given by \cite{cliffe_multilevel_2011}
\begin{equation}\label{eq:OptSamples}
    M_{\ell} = \left\lceil 2\varepsilon^{-2}\sqrt{V_{\ell}/C_{\ell}}\left(\sum_{\ell'=0}^{L}\sqrt{V_{\ell'}C_{\ell'}}\right) \right\rceil\,,
\end{equation}
where $\lceil\cdot\rceil$ is the ceiling function defined as $\lceil a \rceil = \min\{m\in\mathbb{Z}\, |\, m\geq a  \}$ given any $a\in\mathbb{R}$. 

Generating addition samples $\Delta Q_{\ell}^{(i)}$ on level $\ell$ in the line 4 of algorithm~\ref{alg:MLMC}, requires solving uncertain kinetic equation~\eqref{eq:genericKEwU} on levels $\ell$ and $\ell-1$, with the same realization of the random variable $\omega$. Since, in most cases, the variance $V_{\ell}$ is not known a priori, it must be estimated as,
\begin{displaymath}
    V_{\ell} = \mathbb{E}\left[\left\lVert\Delta Q_{\ell} - \mathbb{E}\left[\Delta Q_{\ell}\right] \right\rVert_{H}^{2} \right] \approx \frac{1}{M_{\ell}-1}\left( \sum_{i=1}^{M_{\ell}}\left\lVert \Delta Q_{\ell}^{(i)} \right\rVert^{2}_{H} - \frac{1}{M_{\ell}}\left\lVert\sum_{i=1}^{M_{\ell}} \Delta Q_{\ell}^{(i)} \right\rVert^{2}_{H} \right)\,.
\end{displaymath}
Due to the high dimensionality of the kinetic equations, accumulating samples to estimate the expected value and variance can be infeasible. To avoid this, we use Welford's online algorithm~\cite{welford_note_1962} to compute the expected value and variance on-the-fly on each level.

To test for convergence, in line 7 of algorithm~\ref{alg:MLMC}, we verify if
\begin{displaymath}
    \lVert \mathbb{E}\left[ Q_{L} - Q \right]\rVert_{H} < \varepsilon/\sqrt{2}
\end{displaymath}
is satisfied for the given tolerance $\varepsilon$. Since the exact solution $\mathbb{E}[Q]$ is unknown a priori, we estimate the error based on estimates of the last levels $\Delta Q_{L}$. Assuming that the error of the differences behaves as $\lVert\mathbb{E}[Q_{\ell} - Q_{\ell-1}]\rVert_{H}=\mathcal{O}(2^{-\alpha\ell})$, as $\ell\to\infty$ for $\alpha\geq 1$, we get
\begin{displaymath}
     \lVert \mathbb{E}\left[ Q_{L} - Q \right]\rVert_{H} = \left\lVert \sum_{\ell = L+1}^{\infty}\mathbb{E}\left[ Q_{\ell} - Q_{\ell-1} \right] \right\rVert_{H}\leq \sum_{\ell = L+1}^{\infty}\left\lVert \mathbb{E}\left[ Q_{\ell} - Q_{\ell-1} \right] \right\rVert_{H} \simeq\frac{\left\lVert \mathbb{E}\left[ Q_{L} - Q_{L-1} \right] \right\rVert_{H}}{2^\alpha - 1} \,, \label{eq:conv}
\end{displaymath}
where 
\begin{displaymath}
    A\simeq B \iff cB \leq A \leq \hat{c}B,\quad c,\hat{c}>0\,.
\end{displaymath}

\begin{remark}
    For robustness of the bias estimate, we extrapolate the error from the previous two levels $\left\lVert \mathbb{E}\left[ Q_{L-1} - Q_{L-2} \right] \right\rVert_{H}$ and $\left\lVert \mathbb{E}\left[ Q_{L-2} - Q_{L-3} \right] \right\rVert_{H}$, and check if the maximum of the three satisfies the convergence criterion. Note that since this is a heuristic approach that does not provide theoretical guarantees, we choose a very conservative strategy. Our later convergence tests in Figure~\ref{fig:Accuracy_MC_vs_MLMC} serve motivate this and confirm that it is effective in keeping the error within the prescribed accuracy.
\end{remark}

Note that the RaDLR-MLMC has been constructed to preserve the spatial accuracy of the underlying discretization method. Thus, we can use the following generic MLMC theorem to characterize the total cost of the RaDLR-MLMC estimate $\widehat{Q}_{L}^{\mathrm{ML}}$ in a Hilbert space $H$. 
\begin{theorem}[\cite{giles_multilevel_2015}]\label{thm:MLMC}
    Assume that there exist constants $\alpha,\beta,\gamma >0$, such that $\alpha\geq \frac{1}{2}\min\{\beta,\gamma\}$ and we have the asymptotic rates
    \begin{enumerate}
        \item[(i)] $\lVert \mathbb{E}[Q_{\ell} - Q]\rVert_{H} \leq c_{1}2^{-\alpha\ell} $,
        \item[(ii)] $V_{\ell} \leq c_{2}2^{-\beta\ell} $, and
        \item[(iii)] $C_{\ell} \leq c_{3}2^{\gamma\ell} $.
    \end{enumerate}
    Then for any given $\varepsilon$ small enough, there exists a total number of level $L$ and the number of samples $M_{\ell}$, $\ell = 0,\ldots,L$, such that the MSE satisfies
    \begin{displaymath}
        e\left( \widehat{Q}_{L}^{\mathrm{ML}} \right)^{2} < \varepsilon^{2},
    \end{displaymath}
    and
    \begin{displaymath}
        \mathcal{C}\left( \widehat{Q}_{L}^{\mathrm{ML}} \right) \leq c\begin{cases}
            \varepsilon^{-2}, & \mathrm{if}~\beta>\gamma,\\
            \varepsilon^{-2}(\ln{\varepsilon})^{2}, & \mathrm{if}~\beta=\gamma,\\
            \varepsilon^{-2-\frac{\gamma-\beta}{\alpha}}, & \mathrm{if}~\beta<\gamma.
        \end{cases}
    \end{displaymath}
\end{theorem}
The computational complexity predicted by the MLMC theorem above is guaranteed, provided that the assumptions (i)--(iii) are satisfied. These are often difficult to verify theoretically for complex models, and we refer the reader to the numerical experiments section for numerical verifications for various kinetic equations considered in this paper.

Finally, we note that the classical complexity bound of Theorem~\ref{thm:MLMC} identifies three regimes, depending on the interplay between the variance decay rate $\beta$ and the cost growth rate $\gamma$. The optimal regime $\beta > \gamma$ yields complexity $\mathcal{O}(\varepsilon^{-2})$, matching plain Monte Carlo in the idealized setting where sample cost is fixed, i.e., independent of the discretization level. The borderline case $\beta = \gamma$ incurs only a logarithmic overhead, $\mathcal{O}(\varepsilon^{-2} |\log\varepsilon|^{2})$. While these favorable rates are often attainable through a judicious choice of discretization---for instance, when levels correspond solely to different step sizes in a time-stepping scheme---more involved discretizations combining spatial, temporal, and, as in this work, rank-based refinement often fall into the worst-case regime $\beta < \gamma$, particularly for computationally expensive models. It is important to stress, however, that even in this worst case, MLMC still reduces the complexity exponent by $\beta/\alpha$ relative to plain Monte Carlo, which can translate to savings of several orders of magnitude, as we demonstrate in the numerical experiments section.

\begin{remark}
    In case of multiple QoIs, say $\overline{Q} = (Q_{1},\ldots,Q_{k})$, we estimate the variance and expected difference on level $\ell$ as the maximum of the individual variances and expected differences. That is, 
\begin{displaymath}
\begin{aligned}
    V_{\ell} &:= \max_{i=1,\ldots,k}\left\{ V_{\ell,1},\ldots,V_{\ell,k} \right\}\,, \\
    \left\lVert \mathbb{E}\left[ \overline{Q}_{L} - \overline{Q}_{L-1} \right] \right\rVert_{H} &:= \max_{i=1,\ldots,k}\{\left\lVert \mathbb{E}\left[ Q_{L,i} - Q_{L-1,i} \right] \right\rVert_{H}\}\,,
\end{aligned}
\end{displaymath}
where $V_{\ell,i}$ and $\left\lVert \mathbb{E}\left[ Q_{L,i} - Q_{L-1,i} \right] \right\rVert_{H}$ are the variance and expected differences on level $\ell$ for the $i^{\mathrm{th}}$ QoI. This strategy is a conservative choice that doesn't strictly fit into the setting of Theorem~\ref{thm:MLMC}; however, it preserves the algorithmic guarantees component-wise.
\end{remark}

\section{Numerical results}\label{sec:NumercalExperiments}
In the following, we present numerical experiments ranging from a one-dimensional test case, for which error convergence and computational costs can be assessed against a semi-analytical solution, to a challenging two-dimensional radiation transport benchmark, a high-dimensional radiation therapy dose calculation, and finally the hyperbolic shallow water moment equations, demonstrating the generality of our method beyond radiation transport. In all experiments, the velocity variable is discretized using the $P_{\mathrm{N}}$ method~\cite{case_linear_1967} with $n$ basis functions and the MLMC methods use 10 warm-up samples except the shallow water equations and radiation therapy in which we use 20 and 50, respectively.\footnote{All numerical experiments presented in this section were conducted using our open-source Julia implementation, available at \href{https://github.com/chinsp/publication-RaDLR-MLMC-for-Kinetic-Equations.git}{https://github.com/chinsp/publication-RaDLR-MLMC-for-Kinetic-Equations.git}.}. 

\subsection{Radiation transport}\label{sec:NERTE}
We consider a one-speed radiation transport equation such that $\Omega\in\mathcal{V}=\mathbb{S}^{d-1}$\, is the direction of flight of a particle traveling at speed $v_{0}$. Without loss of generality we assume $v_{0} = 1$. The spatial and angular domains $\mathcal{D}$ and $\mathcal{V}$ and the advection operator $\mathcal{A}(u)\partial_{x}u$ are specified in Table~\ref{tab:RTEDomain} for the one- and two-dimensional slab geometry. Note that in the one-dimensional slab geometry, we replace the angular variable $\Omega$ by $\mu$, which corresponds to the projection of $\mathbb{S}^{2}$ onto $\mathbb{R}$ under symmetry assumptions.
\begin{table}
        \centering
        \begin{tabular}{c|c|c|c}
            Dimensions & $\mathcal{D}$ & $\mathcal{V}$ & $\mathcal{A}(u)\partial_{x}u$ \\
             \hline
            $d=2$ & $[-3,3]^{2}$ & $ \mathbb{P}_{\mathbb{R}^{2}}(\mathbb{S}^{2})$ & $\Omega\cdot\partial_{x}u$\\
            $d=1$ & $[-3,3]$ & $\mathbb{P}_{\mathbb{R}}(\mathbb{S}^{2})$ & $\mu\,\partial_{x}u$ \\
            \hline
        \end{tabular}
        \caption{Domain of the spatial variable $x\in\mathcal{D}$ and the angular variable $\Omega\in\mathcal{V}$ for the radiative transport equation in one-, and two-dimensional slab geometry. The projection of $\mathbb{S}^{2}$ onto $\mathbb{R}^{2}$ is defined as $\mathbb{P}_{\mathbb{R}^{2}}(\mathbb{S}^{2}) \coloneqq \left\{\Omega = (\sqrt{1-\mu^{2}}\sin{\theta},\sqrt{1-\mu^{2}}\cos{\theta})~\Big\vert~\theta\in[0,2\pi),~ 0\leq\mu\leq 1 \right\}$ and $\mathbb{P}_{\mathbb{R}}(\mathbb{S}^{2}) \coloneqq \left\{\mu~\big\vert~ -1\leq\mu\leq 1 \right\}$.}
        \label{tab:RTEDomain}
    \end{table}

The collision operator is assumed to be isotropic. It is given by $ \mathcal{S}(u;\omega) \coloneqq \sigma_{s}\left( \int_{\mathcal{V}} u~\mathrm{d}\Omega - u \right)-\sigma_{a}u\, + s$, where $s$ is the source term. In the collision operator, $\sigma_{s}$ denotes the scattering coefficient and $\sigma_{a}$ the absorption coefficient specifying the probability of scattering and absorption of particles, respectively.

\subsubsection{Gaussian pulse test case}
In the first test case, we consider the uncertain Gaussian pulse test case with an uncertain scattering coefficient~\cite{bennett_uncertainty_2025,bennett_benchmarks_2022,ganapol_analytical_2008}. The QoI for this test case is the scalar flux at time $t = 1$ which we denote by $Q(x) := \int_{-1}^{1} u(t=1.0,x,\mu\,;\omega)\,\mathrm{d}\mu $. The scattering coefficient has the form, $\sigma_{s} = \overline{\sigma}_{s}(1 + 0.1\omega) $, where $\omega\sim\mathrm{U}(-1,1)$ and $\overline{\sigma}_{s}=1.0$ is the scattering coefficient. That is, we assume that there is $10\%$ uncertainty in the probability of scattering of particles. The initial condition is a Gaussian pulse centered at $x=0$
\begin{displaymath}
    u(t=0,x,\mu\,;\omega) = \exp{\left( -x^{2}/\kappa^{2} \right)},
\end{displaymath}
where $\kappa$ is a given constant and $\mathcal{D}$ and $\mathcal{V}$ correspond to the first row in Table~\ref{tab:RTEDomain}. For this problem, a semi-analytic solution is available and is computed using the gPC-sG expansion~\cite{bennett_uncertainty_2025,bennett_benchmarks_2022}. It serves as a reference for computing the errors of the Monte Carlo-based estimators. We denote the semi-analytic solution by $\mathbb{E}[Q_{\mathrm{sa}}]$.

To study the convergence behavior and computational costs of the low-rank based estimators, we compare the RaDLR-MC, full-MLMC (without DLRA), and RaDLR-MLMC with varying tolerances $\varepsilon$ to the semi-analytic solution on a reference grid of size $h_{x,ref} = 5.859375\times10^{-3}$. The coarsest discretization parameters for the MLMC estimators are $h_{x,0} = 0.375$ with $n = 501$ basis functions for angular discretization.

\begin{figure}
    \centering
    \begin{subfigure}[b]{0.45\linewidth}
        \includegraphics[width=\linewidth]{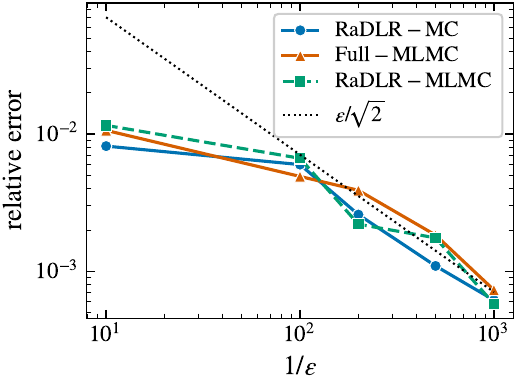}    
    \end{subfigure}
    \begin{subfigure}[b]{0.45\linewidth}
        \includegraphics[width=\linewidth]{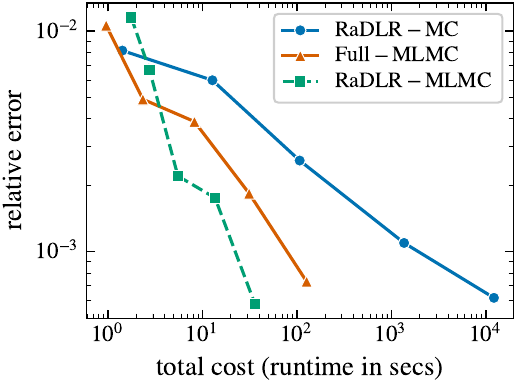}    
    \end{subfigure}
    \caption{Results of the Gaussian pulse test case with $\overline{\sigma}_{s} = 1.0$ and $\kappa = 0.5$ for the RaDLR-MC (with augmented BUG), full-rank MLMC, and RaDLR-MLMC. The solution is compared to the semi-analytic solution $\mathbb{E}[Q_{sa}]$ computed using gPC-sG expansion~\cite{bennett_benchmarks_2022,bennett_uncertainty_2025}. Left: Error of the different estimators with respect to the semi-analytic solution versus the tolerance of the adaptive estimators. Right: Error versus cost (minimum over 5 runs) for the three estimators.}
    \label{fig:Accuracy_MC_vs_MLMC}
\end{figure}

For a fair comparison, we discretize the RaDLR-MC to match the finest level of the full-rank MLMC and the RaDLR-MLMC for a specified tolerance $\varepsilon$. Notably, the finest levels for both estimators coincided across all experiments. The relative error of the estimators and their corresponding costs (minimum over 5 runs) for varying tolerances $\varepsilon$ are shown in Figure~\ref{fig:Accuracy_MC_vs_MLMC}.

From Figure~\ref{fig:Accuracy_MC_vs_MLMC}, we see the bias of all three methods follows the expected curve as $\varepsilon$ is varied. Moreover, by comparing the total computation cost of the estimator we see that the RaDLR-MC, at the finest level, is more expensive than the MLMC estimators at the same relative error. At large tolerances, the full-rank MLMC is computationally cheaper than RaDLR-MLMC which can be attributed to the additional costs of computing QR decompositions in the augmented BUG integrator~\cite{ceruti_rank-adaptive_2022}. At lower tolerances, the finest level of the MLMC estimators is large enough that the computational cost of computing a full-rank sample is much higher than the additional cost of QR decomposition. Thus, the RaDLR-MLMC allows for a computationally cheaper estimate of the expected value of the QoI. These effects are more pronounced in higher-dimensions wherein, for certain problems, it is infeasible to compute a single sample of the full-rank solution.
\subsubsection{Lattice test case}
\begin{figure}[t]
    \centering
    \begin{subfigure}[t]{0.32\linewidth}
    \vspace{0pt}
        \includegraphics[width = 0.9\linewidth]{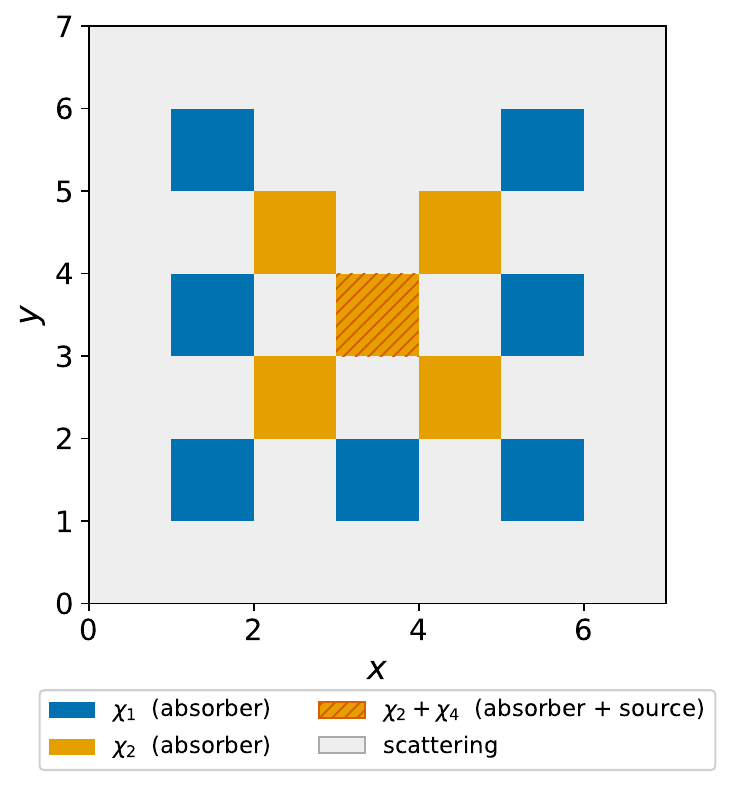}
    \end{subfigure}
    \begin{subfigure}[t]{0.33\linewidth}
        \vspace{0pt}
        \includegraphics[width = \linewidth]{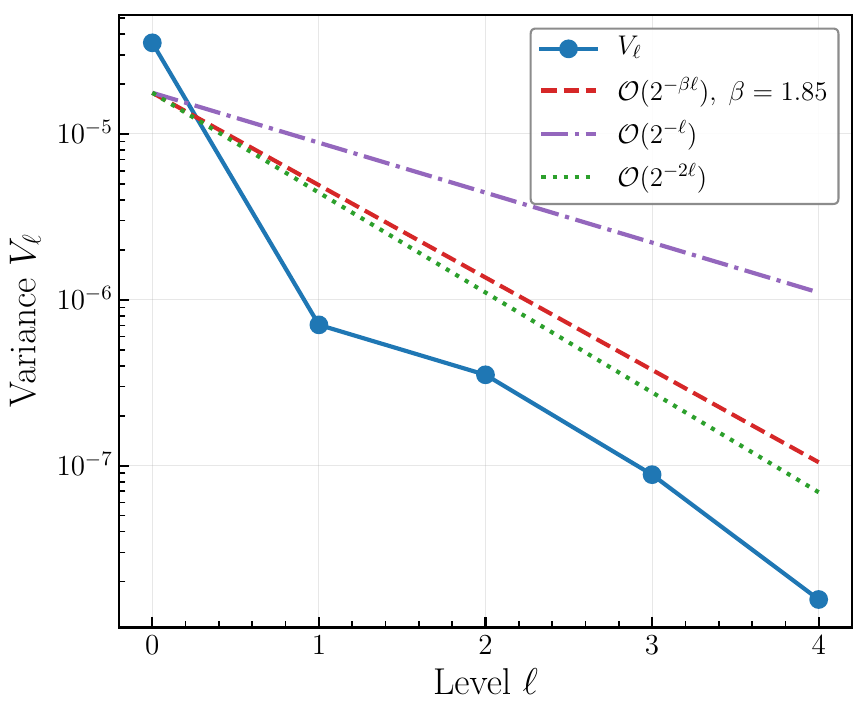}
    \end{subfigure}
    \begin{subfigure}[t]{0.33\linewidth}
    \vspace{0pt}
        \includegraphics[width = \linewidth]{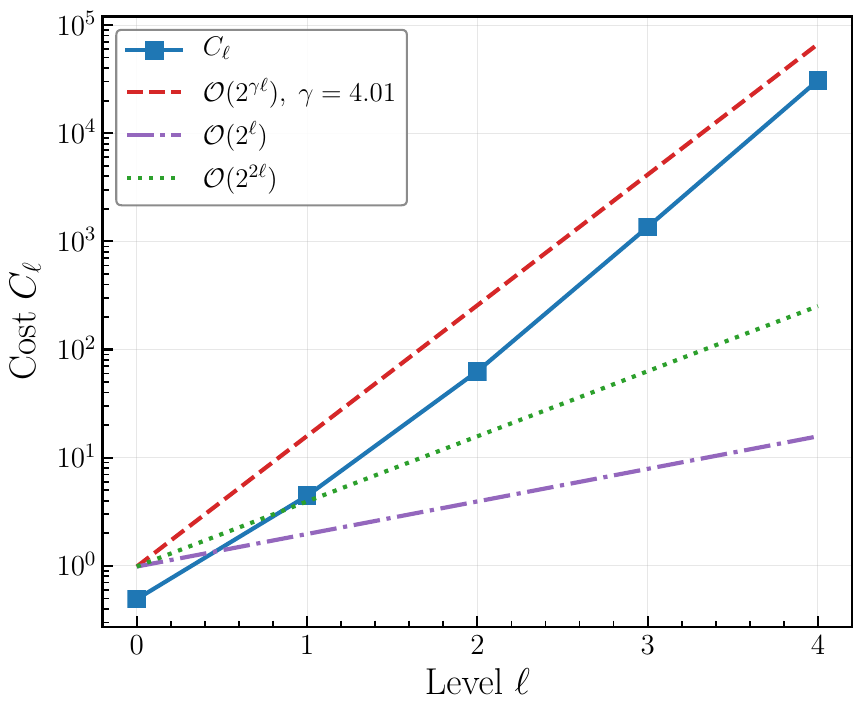}
    \end{subfigure}
    \caption{Left: geometry and setup of the lattice test case with uncertain absorption and source. Middle: variance of the RaDLR-MLMC estimator for the lattice test case across levels. Right: cost of the RaDLR-MLMC estimator for the lattice test case across levels.}
    \label{fig:Lattice_geometry}
\end{figure}

Next, we consider an uncertain lattice test case for the radiation transport equation in two-dimensional slab geometry, which is an important benchmark in the radiative transfer community~\cite{kusch_kit-rt_2023,schotthofer_reference_2025}. The geometry of the lattice is given in Figure \ref{fig:Lattice_geometry}, where the colored blocks represent a strongly absorbing material, while the white regions are treated as weak scatterers of radiation. A constant source placed in the central block emits particles in all directions isotopically. This radiation is absorbed by the blocks, forming radiation shadows away from th source. We provide the details of the random parameters and refer the reader to \cite{schotthofer_reference_2025} for complete details of the test case.

\begin{table}[htbp!]
    \centering
    \begin{tabular}{lccc}
        \hline
        & $\sigma_{a}^{\mathrm{blue}}(\omega)$ 
        & $\sigma_{a}^{\mathrm{orange}}(\omega)$ 
        & $s(t,x,\Omega\,;\omega)$ \\
        \hline
        Experiment 1 
            & $\mathrm{U}(9.5,\;10.5)$ 
            & $\mathrm{U}(9.0,\;11.0)$ 
            & $\mathrm{U}(0.9,\;1.1)$ \\
        Experiment 2 
            & $\mathrm{U}(9.0,\;11.0)$ 
            & $\mathrm{U}(8.5,\;11.5)$ 
            & $\mathrm{U}(0.9,\;1.1)$ \\
        \hline
    \end{tabular}
    \caption{Uncertain parameters for each experiment. The blue outer absorbers, orange inner absorbers, and isotropic source are modelled as independent uniform random variables.}
    \label{tab:experiments}
\end{table}

We assume that the strength of the absorbing blocks in the lattice and the source is uncertain. We divide the absorbing blocks into two regions: the inner absorbers ($\sigma_{a}^{\mathrm{orange}}$; colored in orange in Figure~\ref{fig:Lattice_geometry}) and the outer absorbers ($\sigma_{a}^{\mathrm{blue}}$; colored in blue in Figure~\ref{fig:Lattice_geometry}). The quantity of interest is the scalar flux $Q(x) := \int_{\mathbb{P}(\mathbb{S}^{2})} u(t,x,\Omega\,;\omega)\mathrm{d}\Omega $. The distributions of the random parameters for the test case are given in Table~\ref{tab:experiments}. At the coarsest level, we set $h_{x,0} = 0.21875$ with $1600$ basis functions for angular discretization.

\begin{landscape}
\begin{figure}
    \vspace{-1cm}
    \resizebox{\linewidth}{!}{%
    \begin{tabular}{c@{}ccccc}
    
        \rotatebox{90}{\hspace{3em}$\mathbb{E}[\Delta Q_{\ell}]$} &
        \includegraphics[width=0.19\linewidth]{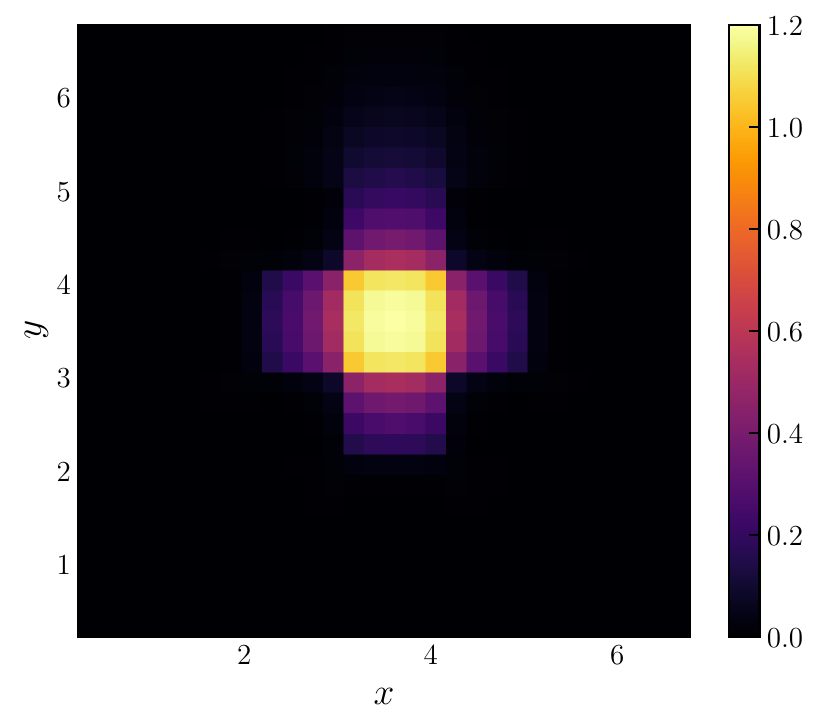} &
        \includegraphics[width=0.19\linewidth]{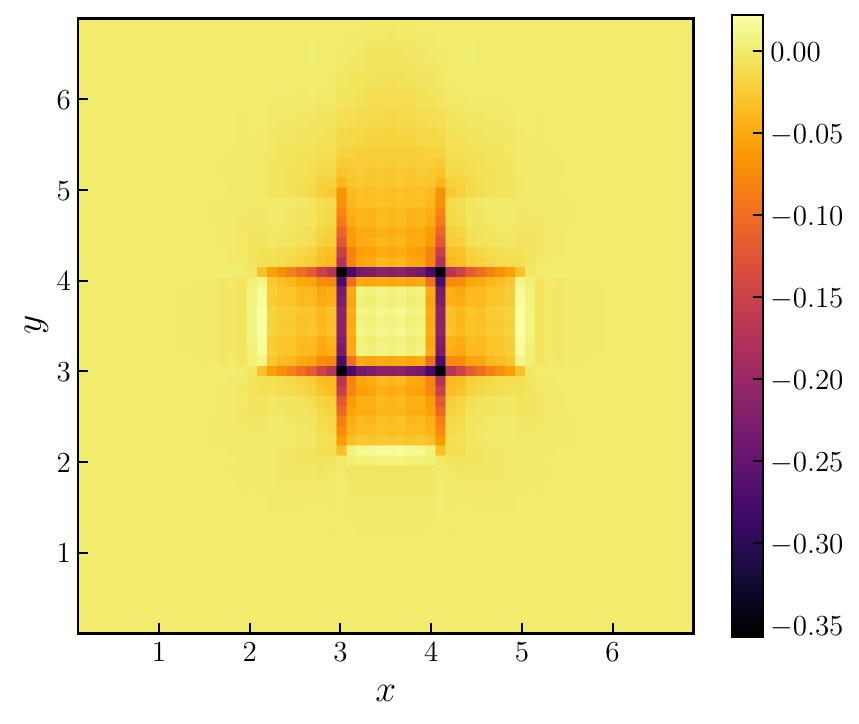} &
        \includegraphics[width=0.19\linewidth]{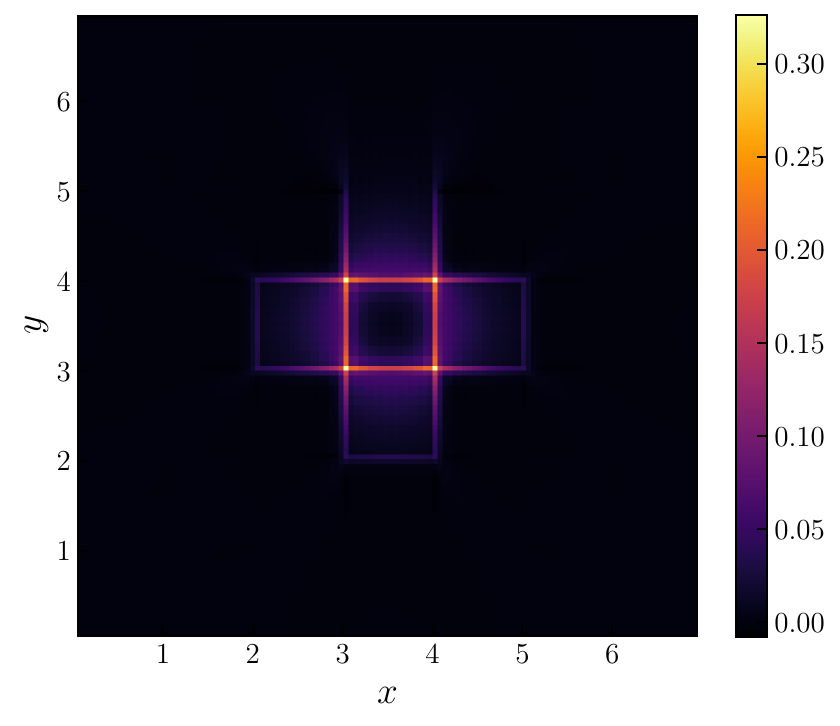} &
        \includegraphics[width=0.19\linewidth]{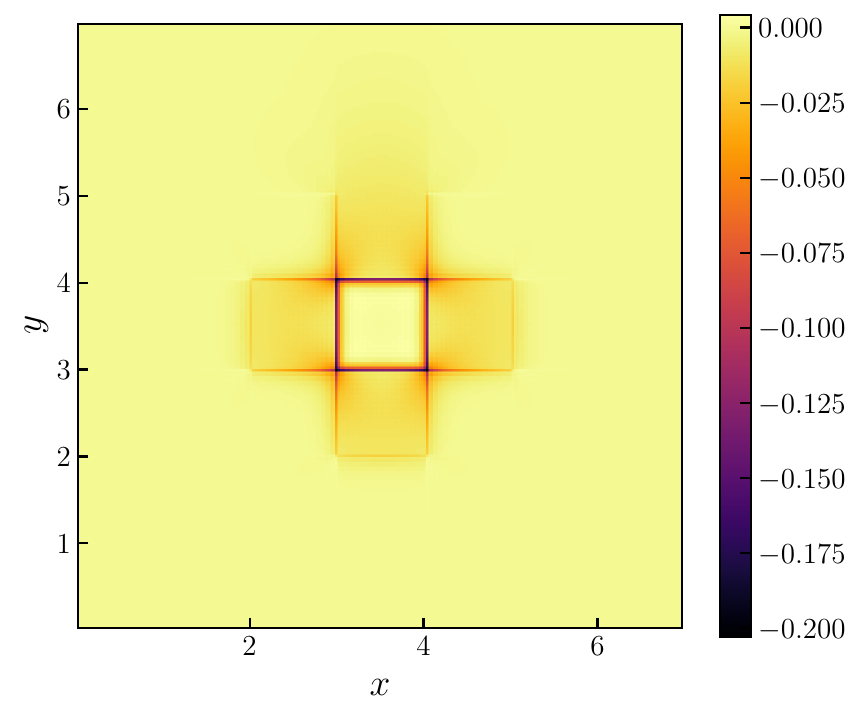} &
        \includegraphics[width=0.19\linewidth]{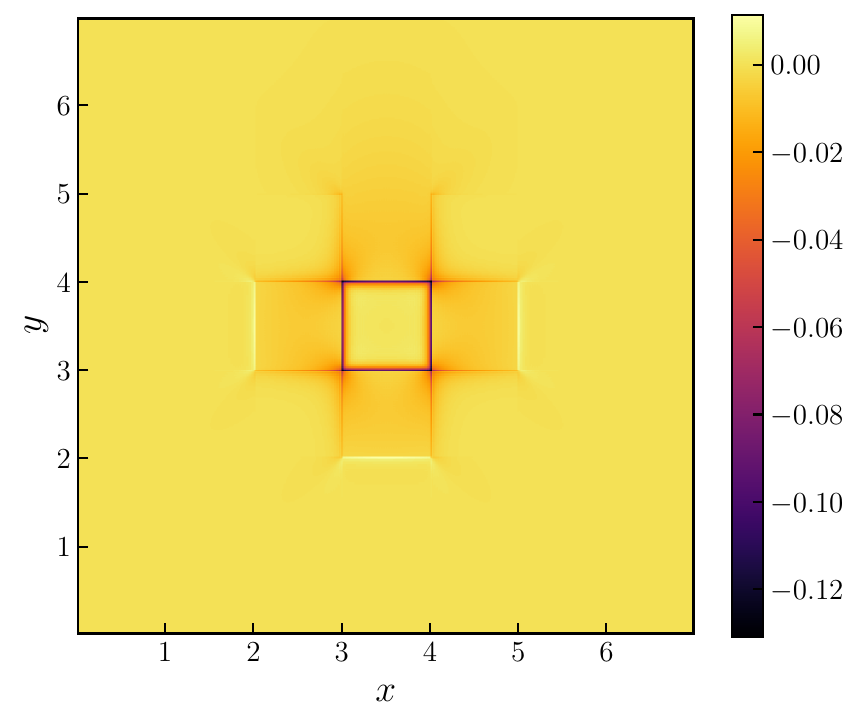} \\
        \rotatebox{90}{\hspace{3em}$\mathbb{E}[Q_{\ell}]$} &
        \includegraphics[width=0.19\linewidth]{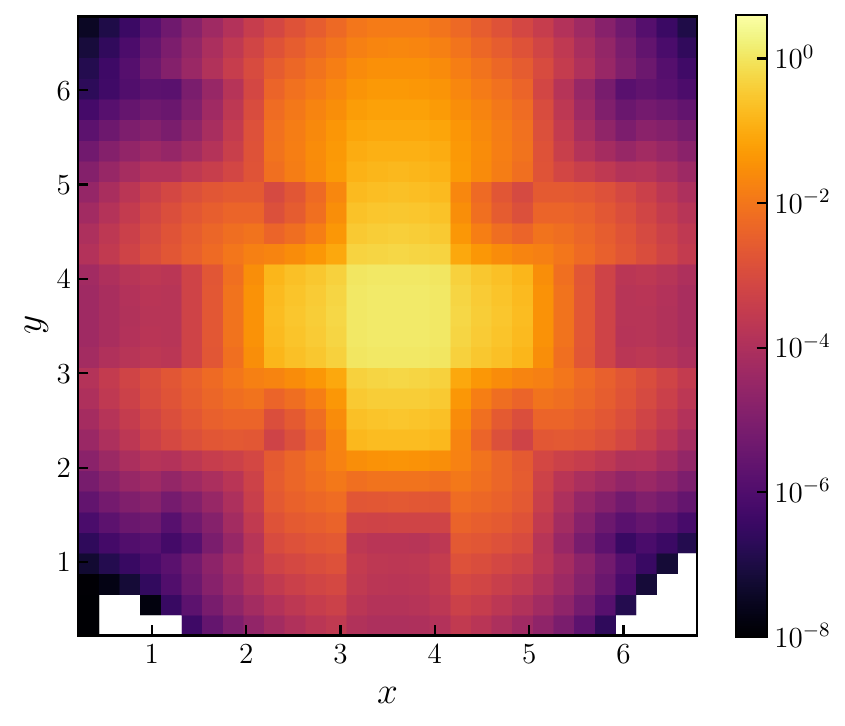} &
        \includegraphics[width=0.19\linewidth]{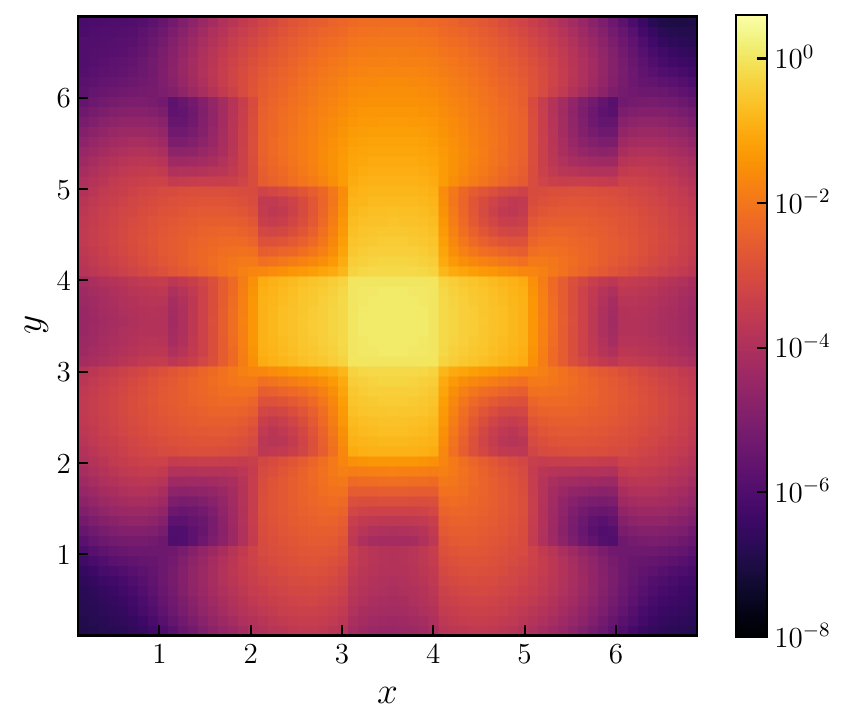} &
        \includegraphics[width=0.19\linewidth]{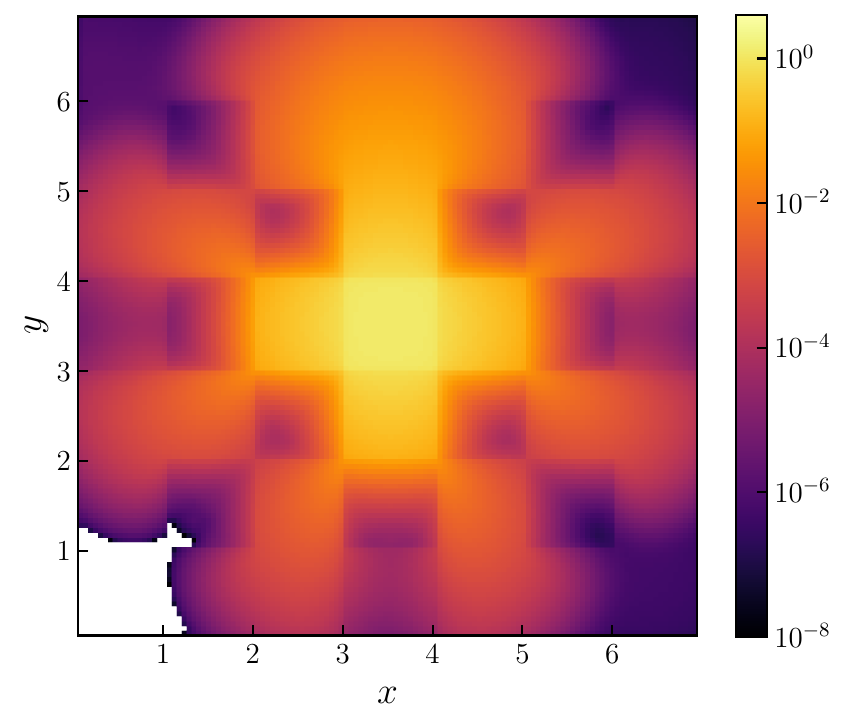} &
        \includegraphics[width=0.19\linewidth]{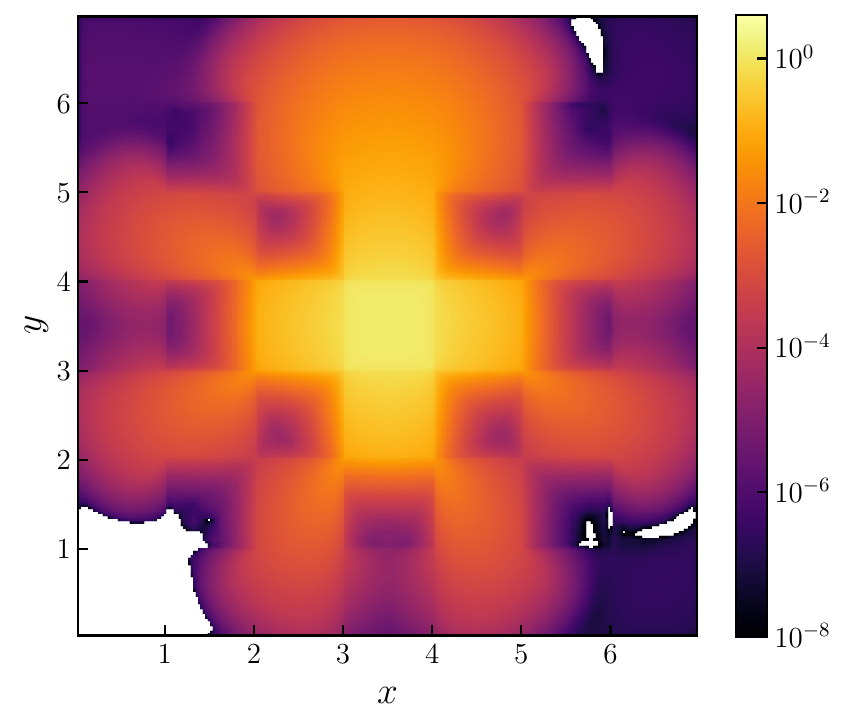} &
        \includegraphics[width=0.19\linewidth]{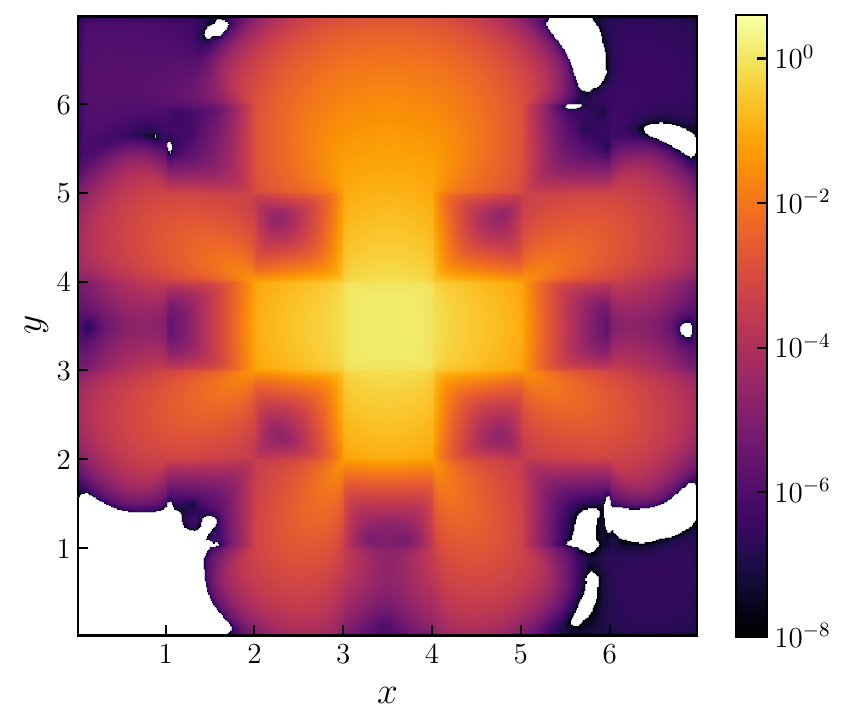} \\
        \rotatebox{90}{\hspace{3em}$\mathbb{E}[\Delta Q_{\ell}]$} &
        \includegraphics[width=0.19\linewidth]{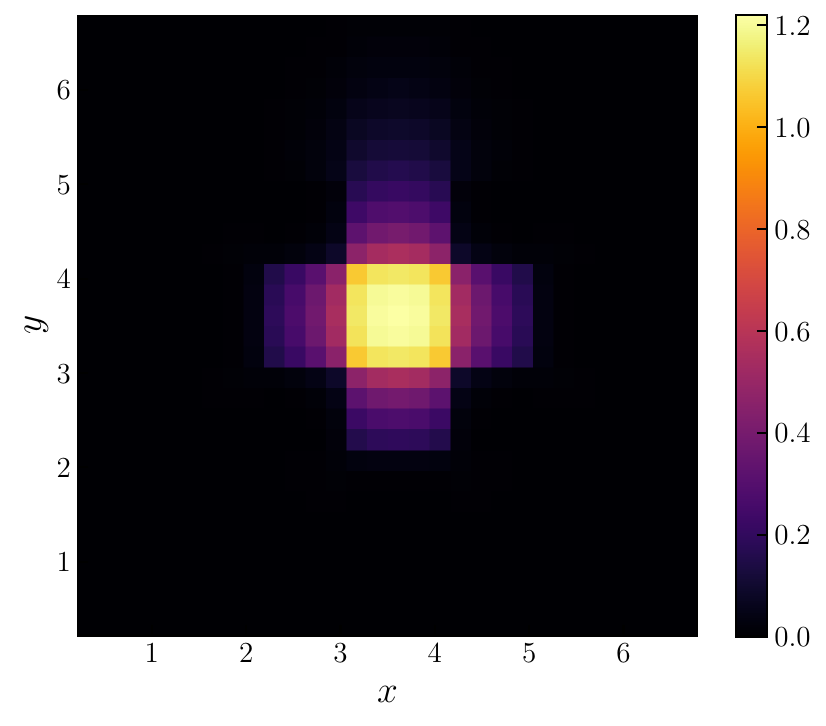} &
        \includegraphics[width=0.19\linewidth]{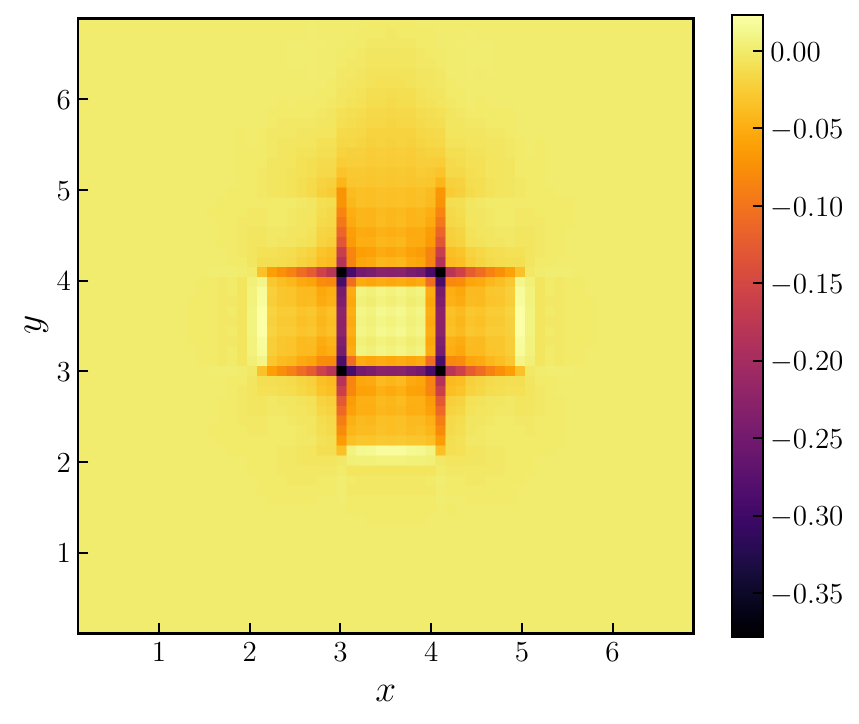} &
        \includegraphics[width=0.19\linewidth]{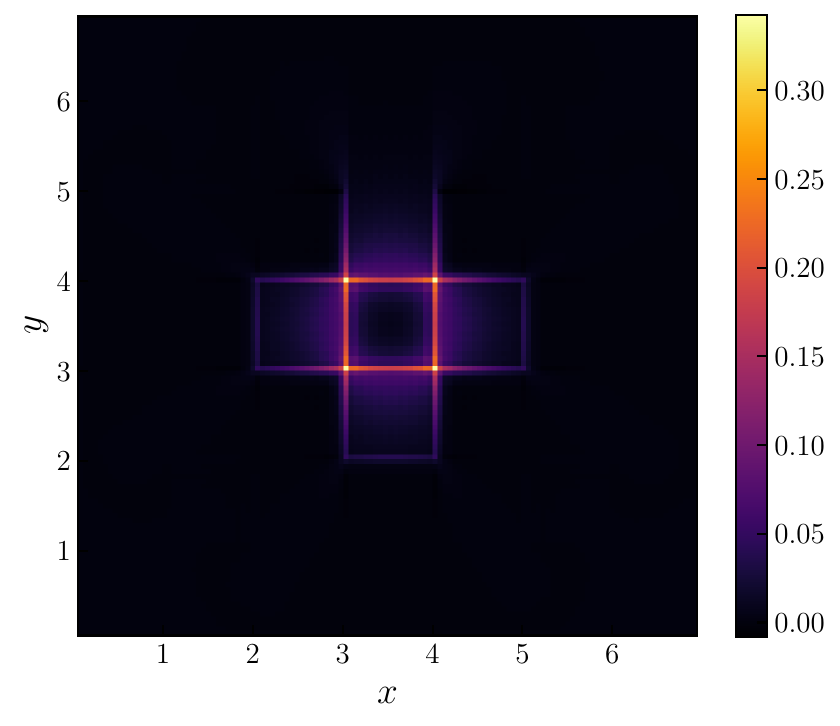} &
        \includegraphics[width=0.19\linewidth]{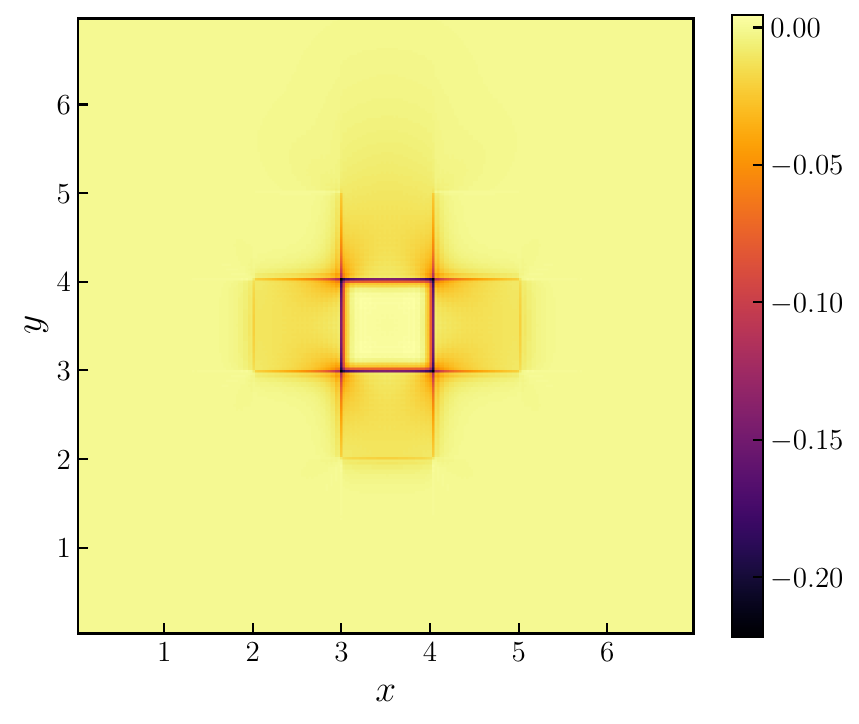} &
        \includegraphics[width=0.19\linewidth]{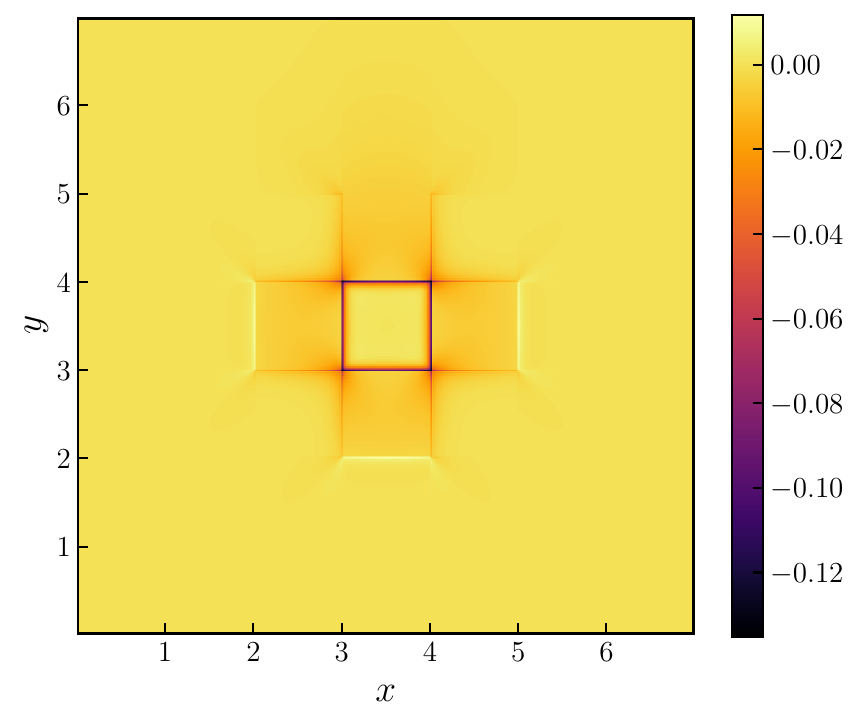} \\
        \rotatebox{90}{\hspace{3em}$\mathbb{E}[Q_{\ell}]$} &
        \includegraphics[width=0.19\linewidth]{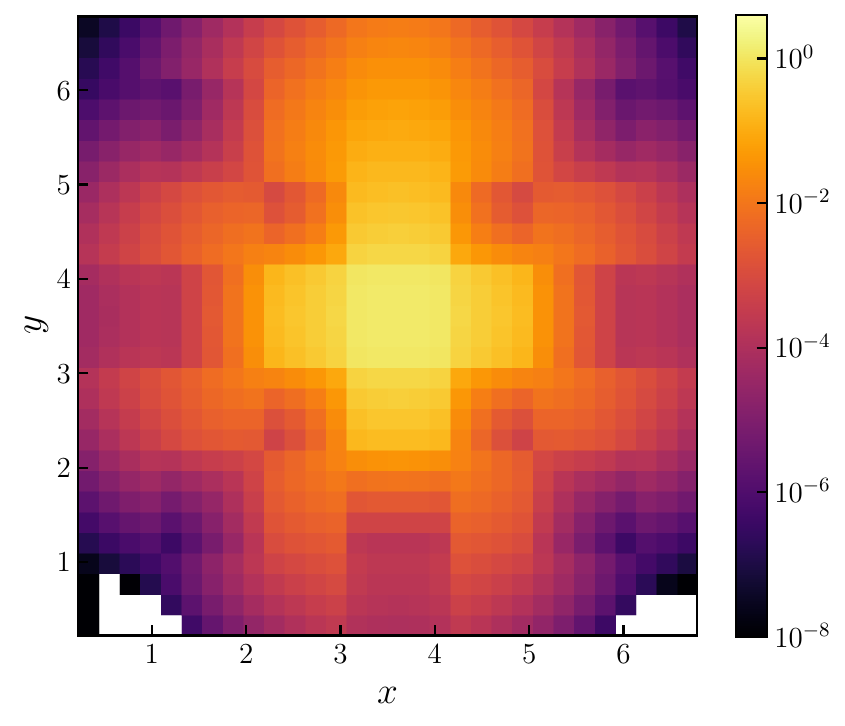} &
        \includegraphics[width=0.19\linewidth]{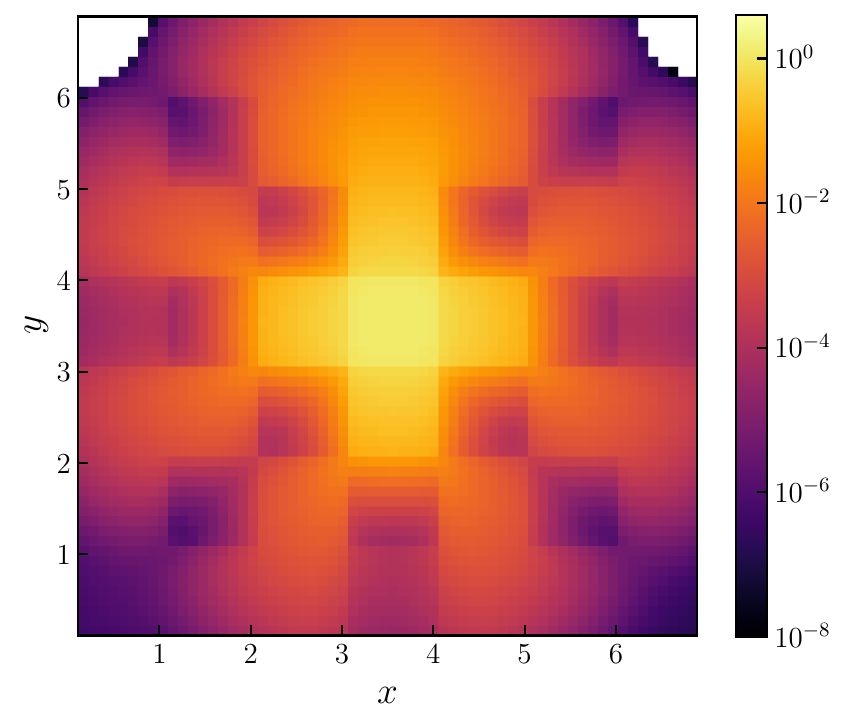} &
        \includegraphics[width=0.19\linewidth]{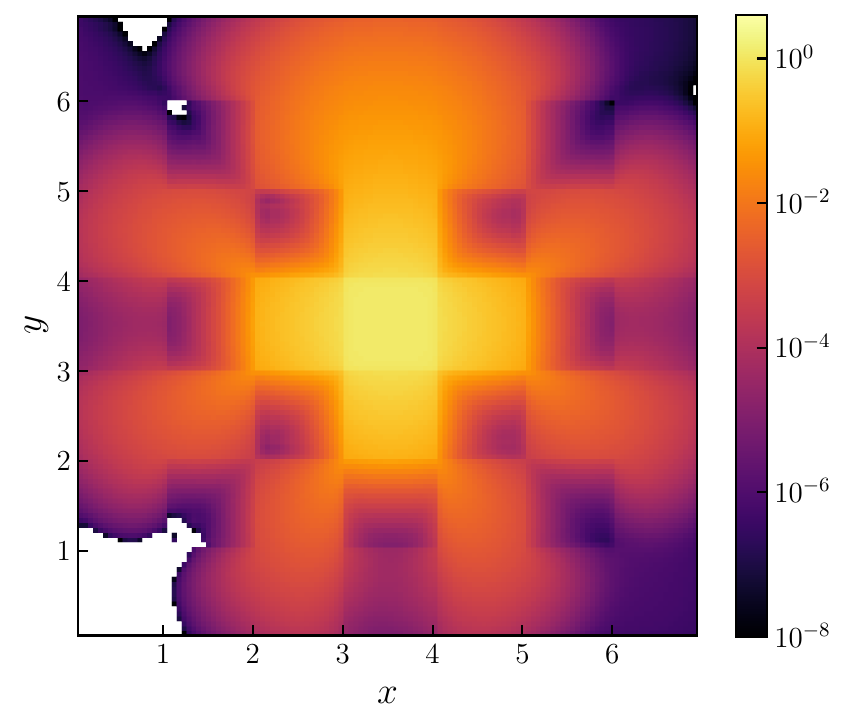} &
        \includegraphics[width=0.19\linewidth]{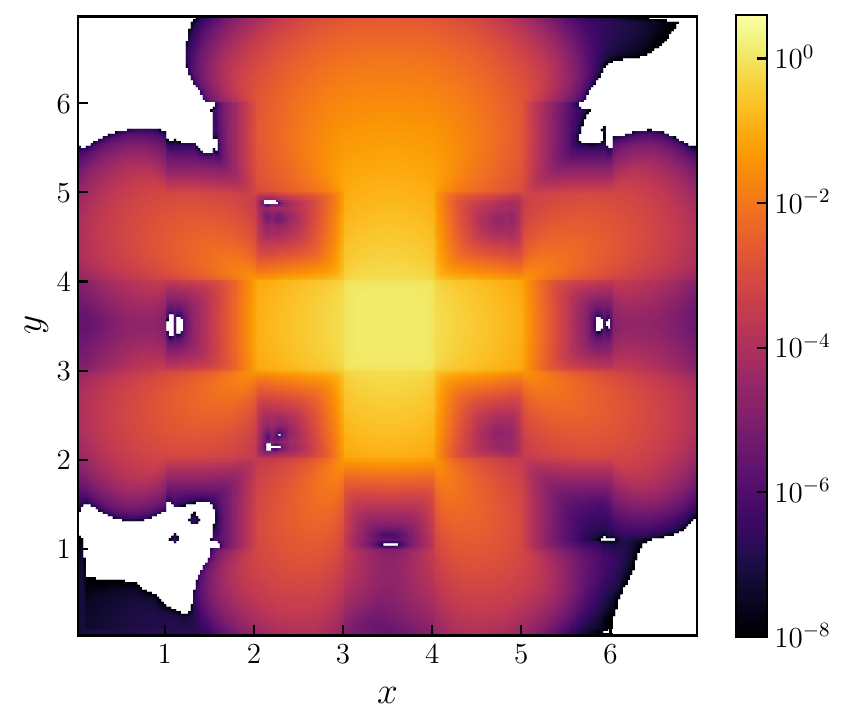} &
        \includegraphics[width=0.19\linewidth]{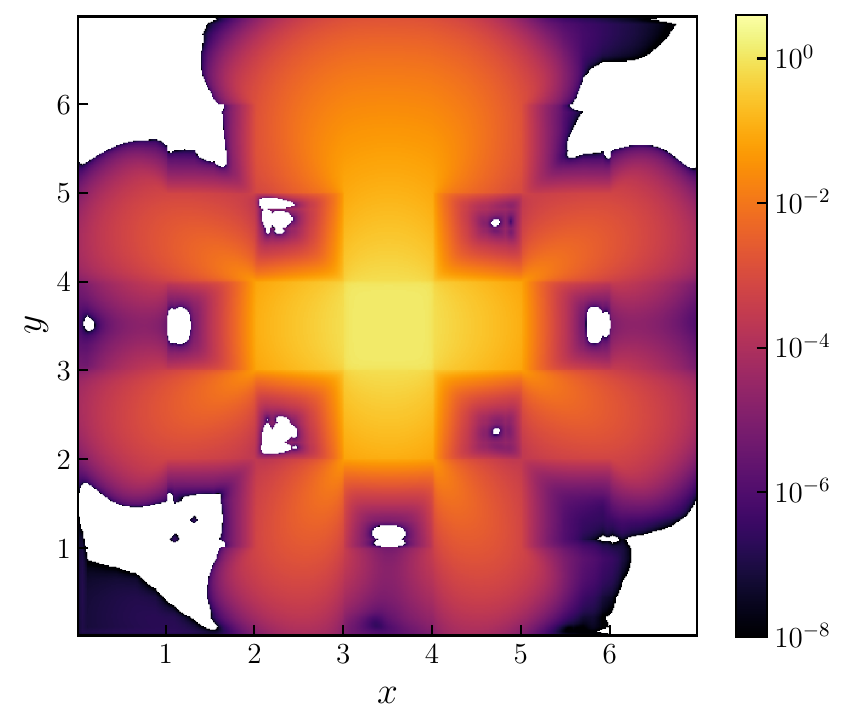} \\
        & $\ell=0$ & $\ell=1$ & $\ell=2$ & $\ell=3$ & $\ell=4$
    \end{tabular}%
    }
    \caption{Quantity of interest (the scalar flux) for the lattice test case with uncertain absorption and source across levels. Row 1,3: $\Delta Q_{\ell}$ for $\ell = 0,1,2,3,4$. Row 2,4: $Q_{\ell} = \sum_{i = 0}^{\ell}\Delta Q_{i}$, $\ell = 0,1,2,3,4$, in a log scale. Row 1,2 correspond to experiment 1 in Table~\ref{tab:experiments} while rows 3,4 are for experiment 2. The white regions correspond to values under the minimum cut-off $10^{-8}$ and negative values.}
    \label{fig:Lattice_results}
\end{figure}
\end{landscape}

In Figure~\ref{fig:Lattice_results}, we plot $\mathbb{E}[Q_{\ell}]$ (in log-scale) and the differences $\mathbb{E}[\Delta Q_{\ell}]$ for the two test cases. We see that the majority of the statistical variation is captured at the coarsest level, while the higher levels correct for the bias in the estimator. We also see that for Experiment 2, with higher uncertainty in the absorbing blocks, the scalar flux shows more oscillations despite being estimated at the same MLMC tolerance $\varepsilon$. In Figure~\ref{fig:Lattice_geometry}, we plot the variance and cost of evaluating a sample on level $\ell$. We see that $\beta<\gamma$, which is not the most favorable case according to the MLMC theorem. However, based on experiments, we see that this drop in efficiency across levels is still not as severe as using the full-rank integrators which have a much higher cost per sample.

\subsection{Radiation therapy}
Next, we consider a more computationally challenging, high-dimensional test case. In radiation therapy, treatment planning requires the computation of energy deposited in a patient's tissue for a specific treatment set up. Patient movements, set up errors, or inaccuracies of available physical and medical data can introduce uncertainties that need to be taken into account \cite{lomax2008intensitya,lomax2008intensityb}. However, already a single simulation of the full six-dimensional transport problem without physical simplifications or model order reduction can be extremely time-consuming. Thus, uncertainty quantification at a high enough resolution is often infeasible in clinical practice \cite{fu2023distributed}.  \\
\begin{figure}[t]
    \centering
    \begin{subfigure}[b]{0.35\linewidth}
        \includegraphics[width = \linewidth]{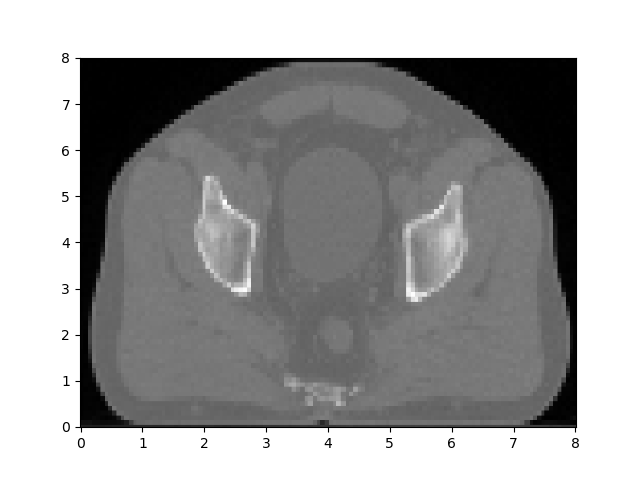}
    \end{subfigure}
    \begin{subfigure}[b]{0.31\linewidth}
        \includegraphics[width = \linewidth]{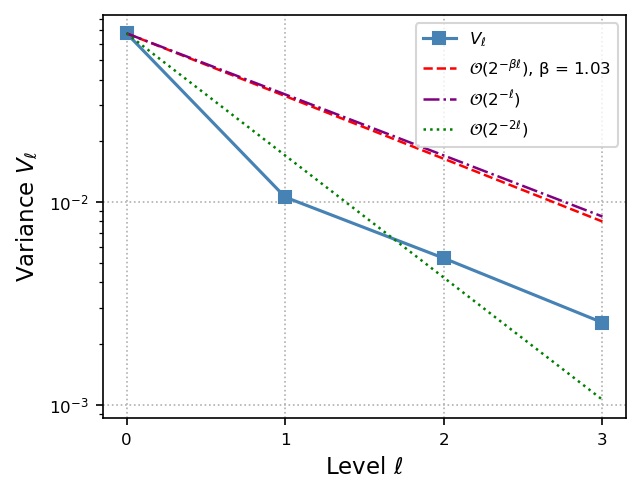}
    \end{subfigure}
    \begin{subfigure}[b]{0.31\linewidth}
        \includegraphics[width = \linewidth]{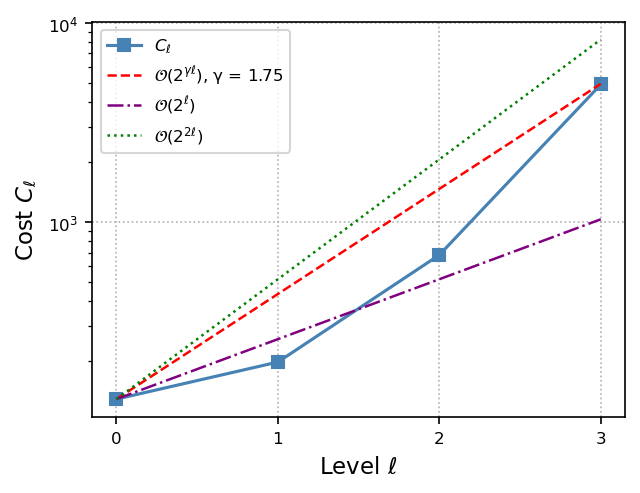}
    \end{subfigure}
    \caption{Left: CT image. Middle: variance of the RaDLR-MLMC estimator for the radiation therapy test case across levels. Right: cost of the RaDLR-MLMC estimator for the radiation therapy test case across levels.}
    \label{fig:RadTherapy_geometry}
\end{figure}
 We now consider a simplified prostate patient CT from the CORT data set \cite{craft2014shared} (see Figure \ref{fig:RadTherapy_geometry}), irradiated with two 75 MeV proton beams which are placed in opposing lateral positions, i.e., at 90° and 270° gantry angle. Further, we assume a positioning error following a normal distribution with standard deviation 3mm. Note, that this choice of irradiation angles and uncertainty is similar to what would be used in application \cite{perko2016fast,cao2015improved}, however we have scaled down the size of the CT scan and adapted the energy of the beams accordingly to limit computational costs. 
 
 Mathematically, we model proton transport using the continuous slowing down approximation to the steady linear Boltzmann transport equation: 
\begin{align}
    \mathbf{\Omega}\cdot\nabla_{\mathbf{r}}\psi(E,\mathbf{r},\mathbf{\Omega}) &- \frac{\partial \mathcal{S}(E,\mathbf{r}) \psi(E,\mathbf{r},\mathbf{\Omega})}{\partial E} - \frac{1}{2}\frac{\partial^2 T(E,\mathbf{r})\psi(E,\mathbf{r},\mathbf{\Omega})}{\partial E^2}= \Gamma \psi(E,\mathbf{r},\mathbf{\Omega}),
\end{align}
where the phase space of the particle density $\psi$ consists of energy $E\in[E_{\mathrm{min}},E_{\mathrm{max}}]\subset \mathbb{R}_{\geq0}$, space $\mathbf{r}\in \mathbb{R}^3$ and direction of flight $\mathbf{\Omega} \in \mathbb{S}^2$.
The stopping power $\mathcal{S}$ describes the average rate of energy loss and straggling $T$ represents the statistical spread in energy loss due to stochastic interactions. The collision operator $\Gamma$ on the right-hand side depends on the differential and total scattering cross sections $\Sigma_s$ and $\Sigma_t$ and is defined as
\begin{align}
    \Gamma \psi = \int_{\mathbb{S}^2} \Sigma_s(E,\mathbf{r},\mathbf{\Omega}'\cdot \mathbf{\Omega})\psi(E,\mathbf{r},\mathbf{\Omega}')\mathrm{d}\mathbf{\Omega}' - \Sigma_t(E,\mathbf{r})\psi(E,\mathbf{r},\mathbf{\Omega}).
\end{align}
  We model the incoming physical beams with a standard deviation of 0.3 cm in space and 1\% in energy. The random lateral beam shifts $\boldsymbol{\omega} \in \mathbb{R}^2$ are applied in beam's eye view (perpendicular to the beam directions) and follow a normal distribution $\mathcal{N}(\boldsymbol{\omega};\boldsymbol{\mu}_{\omega},\boldsymbol{C}_{\omega} )$ with parameters
  \begin{align*}
      \boldsymbol{\mu}_{\omega}= (0,0)^T, \boldsymbol{C}_{\omega} = \begin{bmatrix} 0.3^2 & 0 \\ 0 & 0.3^2 \end{bmatrix}.
  \end{align*} The resolution on the coarsest level is chosen as 0.2 cm in beam direction (along the y-axis) and 0.4 cm in the other dimensions. More details on the implementation and numerical methods used can be found in \cite{stammer2026high}. 
  
\begin{figure}
    \centering
    \begin{tabular}{c@{}cccc}
    
        \rotatebox{90}{\hspace{2em}$\mathbb{E}[\Delta Q_{\ell}]$} &
        \includegraphics[width=0.22\linewidth]{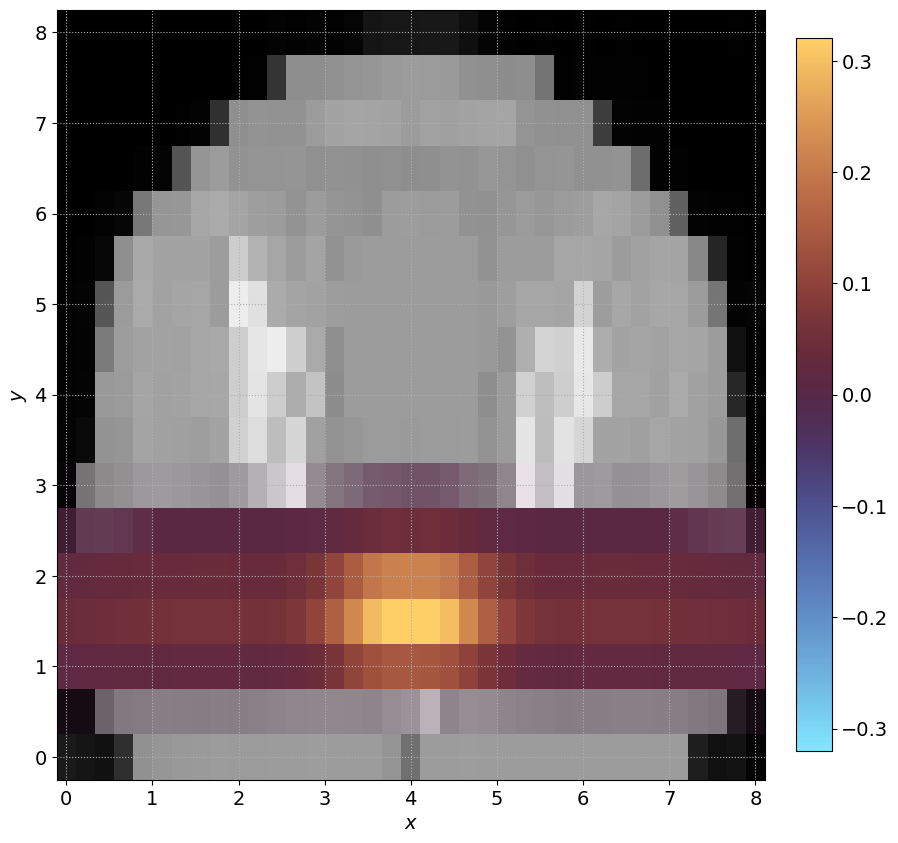} &
        \includegraphics[width=0.22\linewidth]{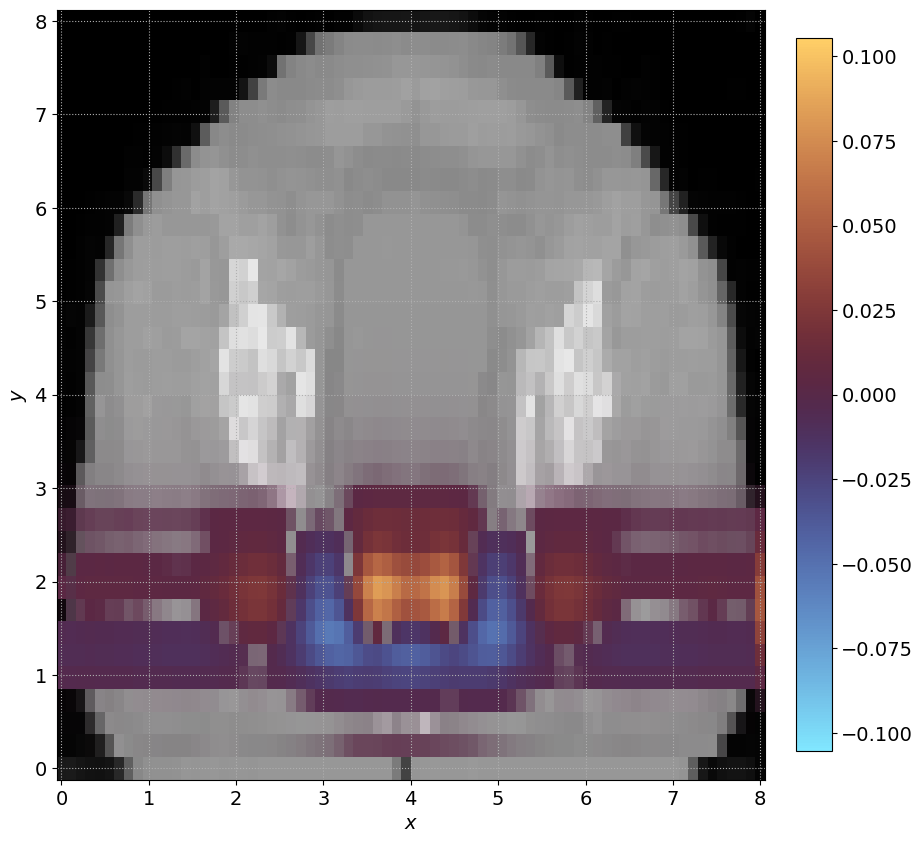} &
        \includegraphics[width=0.22\linewidth]{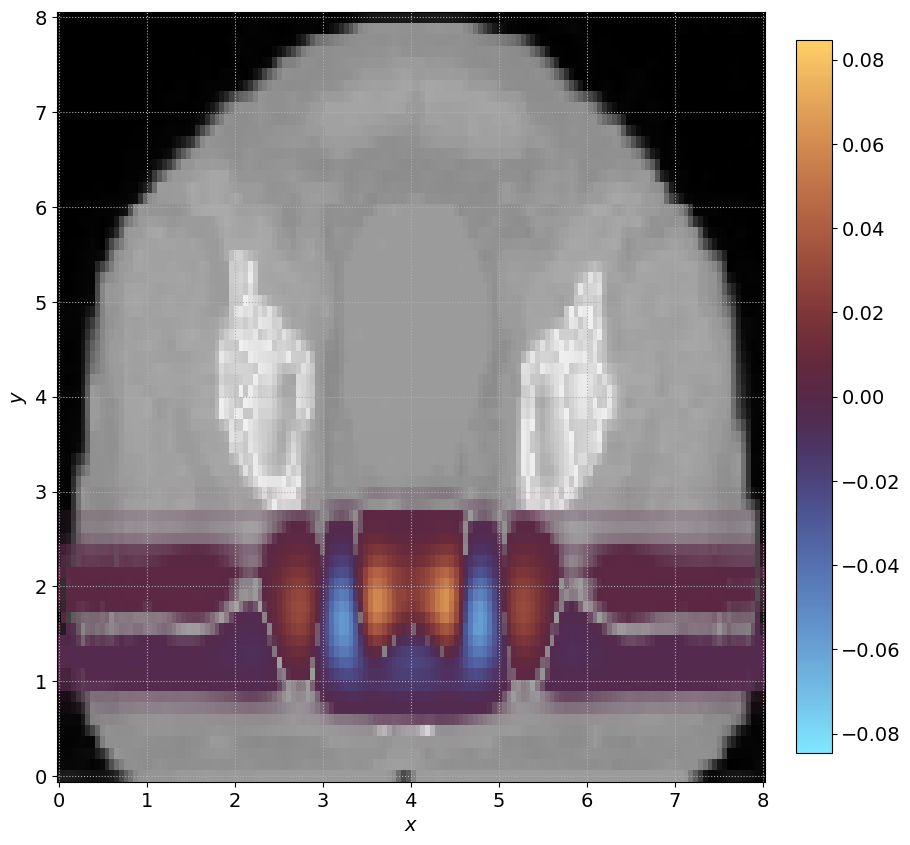} &
        \includegraphics[width=0.22\linewidth]{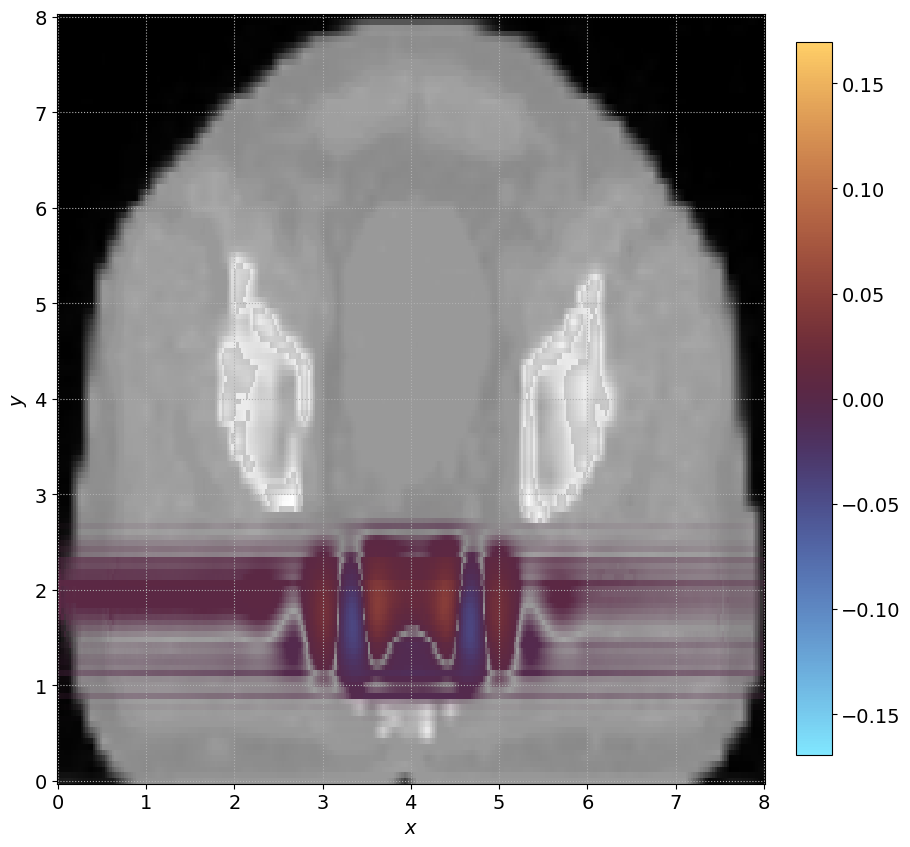} \\
        \rotatebox{90}{\hspace{2em}$\mathbb{E}[Q_{\ell}]$} &
        \includegraphics[width=0.22\linewidth]{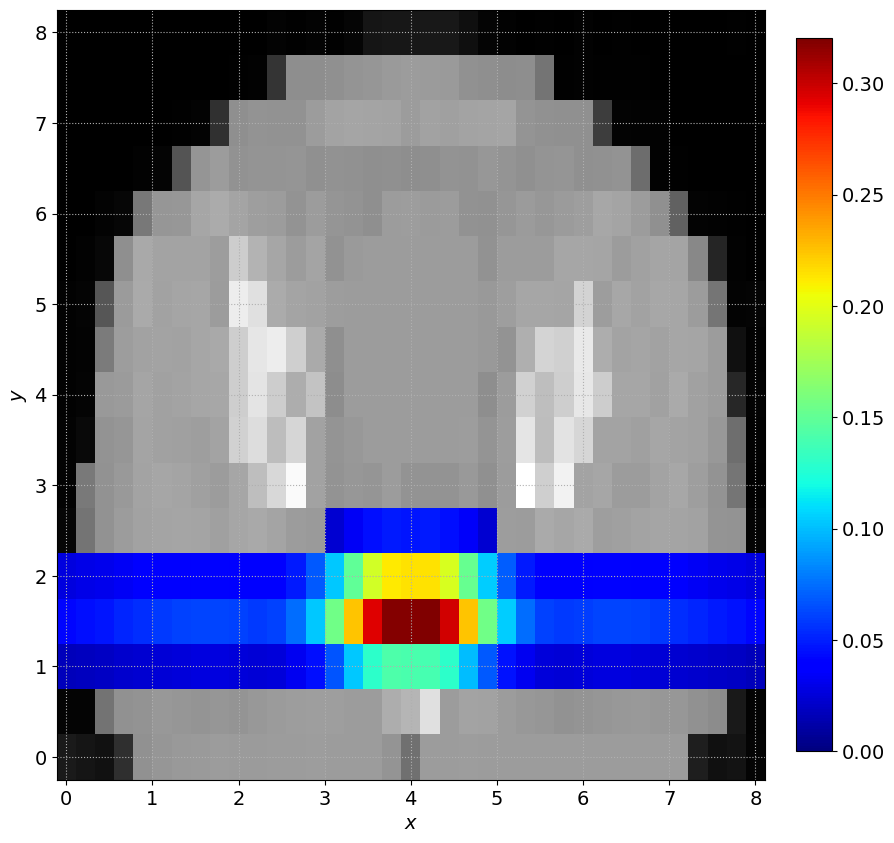} &
        \includegraphics[width=0.22\linewidth]{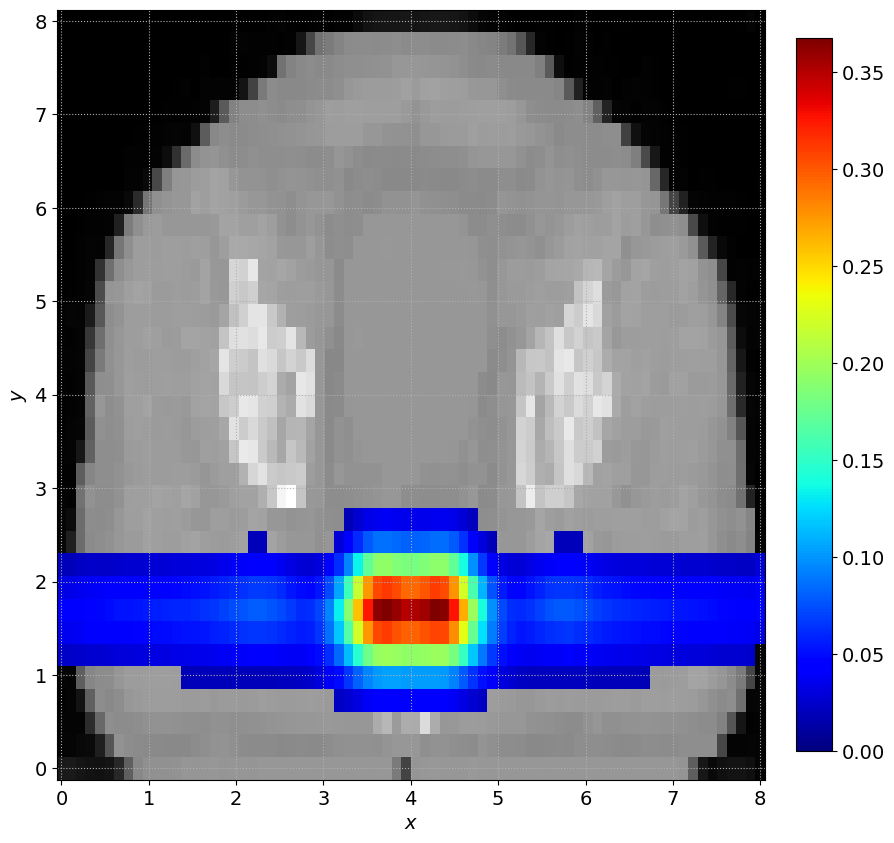} &
        \includegraphics[width=0.22\linewidth]{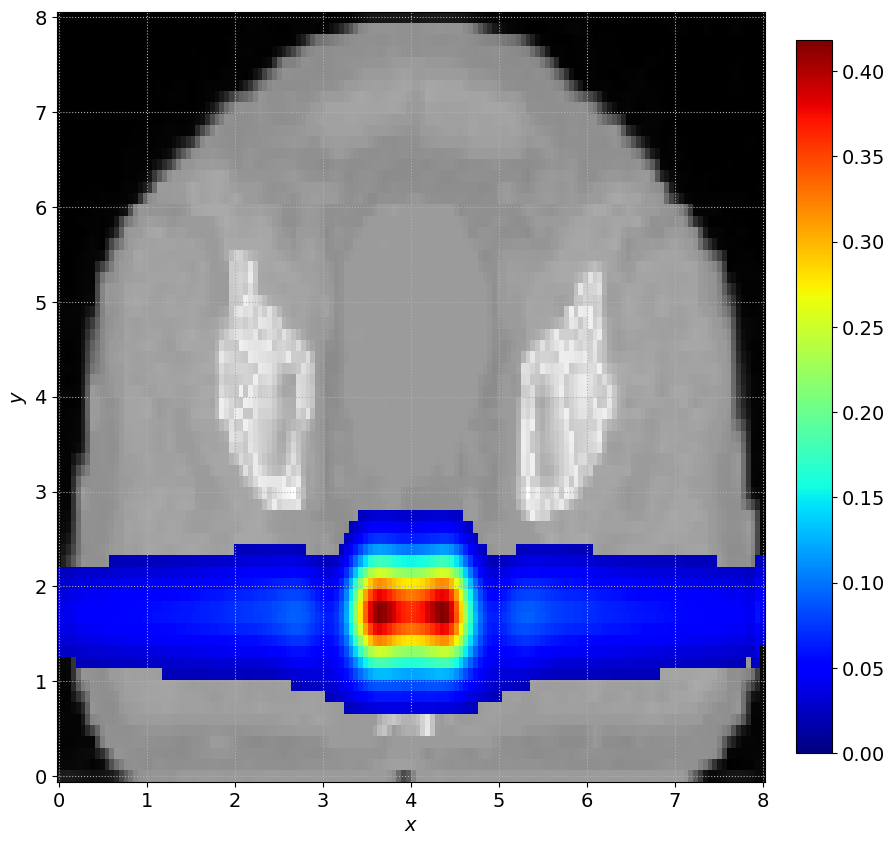} &
        \includegraphics[width=0.22\linewidth]{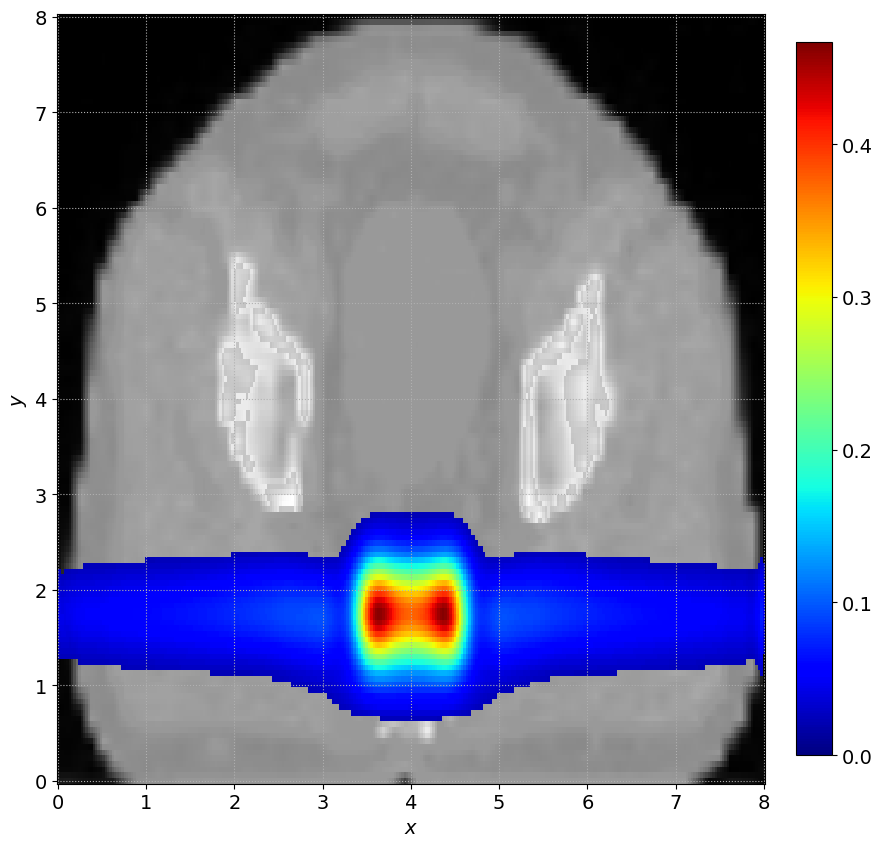} \\

        & $\ell=0$ & $\ell=1$ & $\ell=2$ & $\ell=3$ 
    \end{tabular}

    \caption{Two-dimensional slice along x-y-plane of the quantity of interest (deposited energy) for the radiation therapy test case with uncertain beam positions. Row 1: $\Delta Q_{\ell}$ for $\ell = 0,1,2,3$. Row 2: $Q_{\ell} = \sum_{i = 0}^{\ell}\Delta Q_{i}$, $\ell = 0,1,2,3$. }
    \label{fig:RadTherapy_results}
\end{figure}

Figure \ref{fig:RadTherapy_geometry} shows the costs per sample and the variance for each level. Similarly to the previous test case here $\beta = 1.03 <  1.75 = \gamma$ which theoretically places us in the least favorable case of the MLMC theorem. However, it is clear from Figure \ref{fig:RadTherapy_geometry} that computing a large number of samples for a Monte Carlo estimate of the expected value at highest resolution is very expensive despite already having greatly reduced the costs by use of the dynamical low-rank approximation. Using the low-rank MLMC approach, only a very small number of samples is required at highest resolution, thus making expected value computations much more feasible. Figure \ref{fig:RadTherapy_results} shows the expected energy delivered in the described CT test case as well as the differences and means at each level. Here we can clearly see that while the basic structure of the solution can be represented at coarser levels, the finest levels are necessary to capture essential details in the solution such as the distinct peaks of the two beams. 
\subsection{Hyperbolic Shallow Water Moment Equations}
The hyperbolic shallow water moment equations (HSWME)~\cite{kowalski_moment_2019,koellermeier2020analysis} model free-surface flows arising in applications such as flood prediction, tsunami simulation, and weather forecasting. The HSWME are derived from the Navier--Stokes equations when assuming a polynomial velocity profile along the water depth \cite{kowalski_moment_2019}. The resulting kinetic equations are then regularized to guarantee hyperbolicity as proposed in \cite{koellermeier2020analysis}. Despite being cheaper than simulating incompressible Navier-Stokes equations, the HSWME yet remains expensive to simulate at scale, motivating the use of model order reduction techniques such as DLRA. For a complete model description and details on DLRA for HSWME, we refer the reader to \cite{koellermeier_macro-micro_2024}.

The macroscopic variables in HMSWE are, the height of the water column $h(t,x\,;\omega)$, the velocity profile $p(t,x,z)$, and the momentum of the water flow $hp_{m}$, where $p_{m}(t,x\,;\omega)$ denotes the mean velocity of the water flow. We consider an uncertain dam break test case~\cite{koellermeier_macro-micro_2024} with independent uncertainty in the initial background height, shock position, and shock amplitude and set the QoI to be $Q:=(h,hp_{m})^{\top}$. To e specific, the initial condition is given by
\begin{displaymath}
    h(t=0,x\,;\omega) = 0.3\omega_{1}+ 0.35\omega_{2}\cdot(\tanh{(50x)} - \tanh{(50(x-0.2\,\omega_{3}))})\,,
\end{displaymath}
where $\omega_{i}\in\mathrm{U}(1-\delta/100,1+\delta/100)$, for $i = 1,2,3$, represents $\pm\delta\%$ parametric uncertainty. The flow is assumed to be initially at rest, i.e., $p(t=0,x,z) = 0$. 

\begin{figure}
    \centering
    \begin{subfigure}[b]{\linewidth}
        \includegraphics[width=\linewidth]{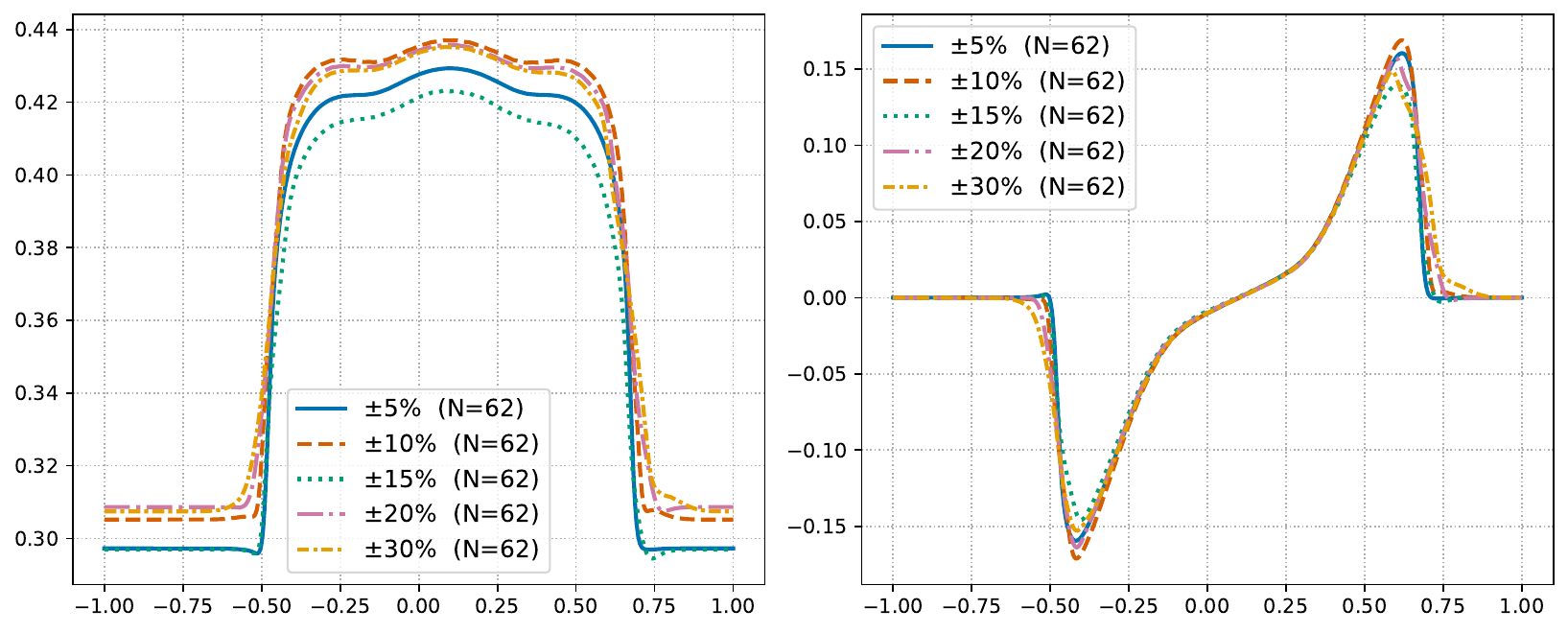}
    \end{subfigure}
    \begin{subfigure}[b]{\linewidth}
        \includegraphics[width=\linewidth]{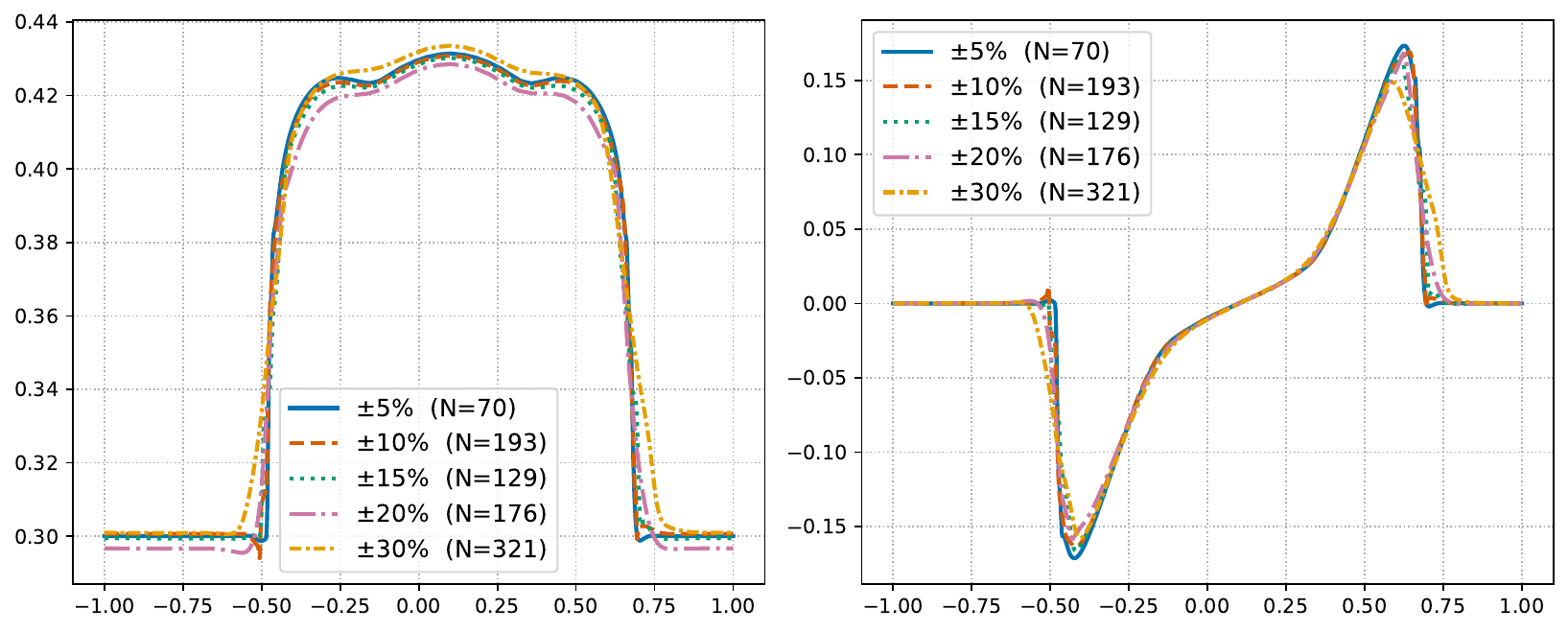}
    \end{subfigure}
    \begin{subfigure}[b]{\linewidth}
        \includegraphics[width=\linewidth]{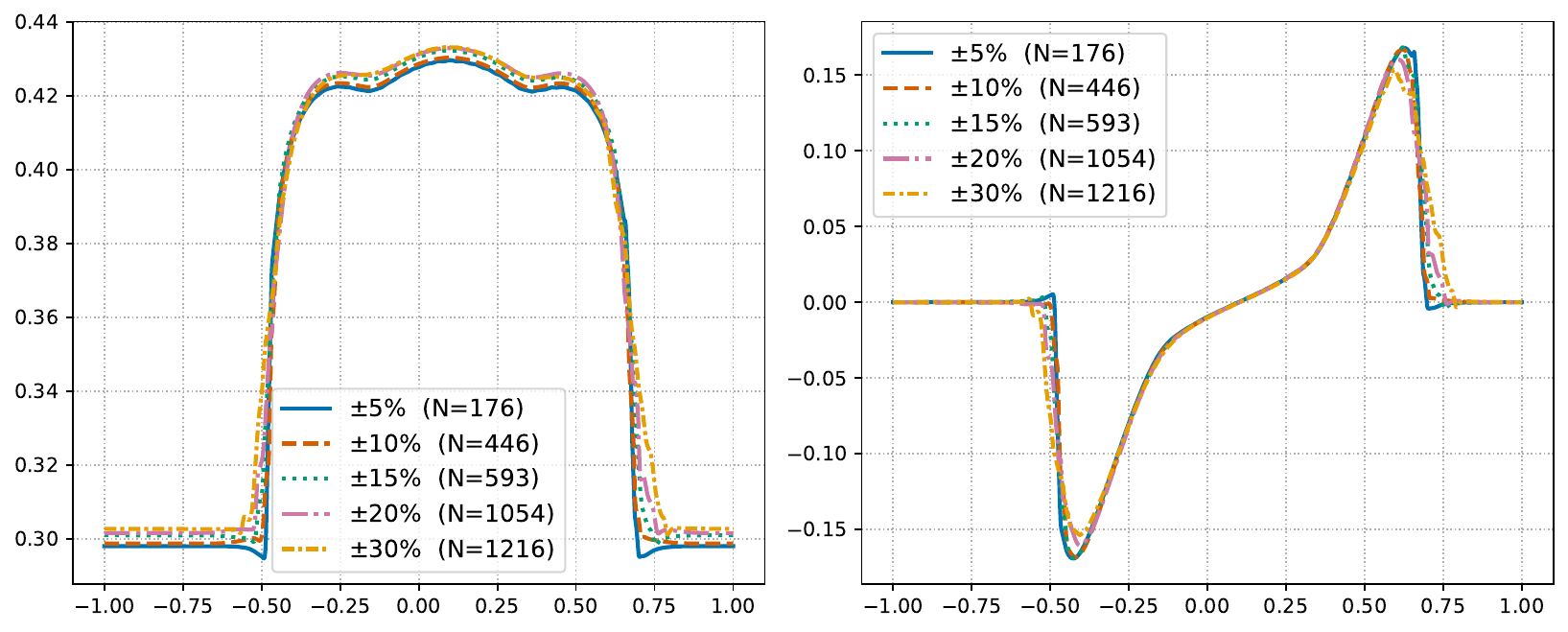}
        \caption*{$\mathbb{E}[h] \qquad\qquad\qquad\qquad\qquad\qquad\qquad\qquad\qquad \mathbb{E}[hu]$}
    \end{subfigure}
    \caption{MLMC estimates of the water height $\mathbb{E}[h]$ (left) and momentum $\mathbb{E}[hp_{m}]$ (right) for the shallow water shock problem. Each curve corresponds to a different level of parametric uncertainty $\pm \delta\%$ of the nominal values, ranging from $\pm 5\%$ to $\pm 30\%$ ($\pm \delta\%$). Results are shown for MLMC tolerances $\varepsilon = 0.1$ (top), $\varepsilon = 0.01$ (middle), and $\varepsilon = 0.005$ (bottom). The number of solver calls $N$ required to meet each tolerance is reported in the legend.} 
    \label{fig:SWE_uncertainy_study}
\end{figure}

The coarsest level is discretized by $h_{x,0} \approx 0.02$ with $102$ basis functions for the angular discretization using the method of moments~\cite{case_linear_1967}. In Figure~\ref{fig:SWE_uncertainy_study}, we plot the expected value of the height and momentum profiles for the dam break test case for three different MLMC tolerances $\varepsilon=0.1$, $\varepsilon = 0.01$, and $\varepsilon = 0.005$. These plots show that the uncertainty in the initial condition of the dam break test case significantly effects the expected outcome of the model. These can be simulated and quantified efficiently with the RaDLR-MLMC estimators by choosing an appropriate tolerance. As expected, reducing $\varepsilon$ improves the estimate of $\mathbb{E}[Q]$ while only minor effects are seen in the expected velocity profile $\mathbb{E}[hp_{m}]$.
\section{Conclusions and outlook}
To overcome the computational bottleneck of sampling-based UQ methods for high-dimensional kinetic equations, we introduced a rank-adaptive dynamical low-rank Monte Carlo estimator that computes a low-rank approximation of the quantity of interest for each sample. A detailed bias analysis revealed a key computational challenge inherent to rank-adaptive integrators: the saturation of errors at smaller tolerances. To address this challenge, we developed a spatially accurate rank-adaptive integrator and leveraged it to construct a rank-adaptive dynamical low-rank multilevel Monte Carlo estimator, significantly reducing the overall computational cost.
We demonstrated the efficacy of the proposed method on a range of challenging test cases. For higher-dimensional problems, the ability to obtain reliable estimates of the quantities of interest at all represents a significant achievement although the variance decay and cost growth remain far from ideal. To facilitate reproducibility and broader adoption, we provide an open-source Julia implementation of the method, compatible with any numerical solver, with or without DLRA.



\section*{Acknowledgements}
The authors would like to thank Martin Frank for his helpful comments and suggestions.

\bibliographystyle{plain}

\bibliography{SMAI-JCM-Online-Template/Final_manuscript}

\begin{thebibliography}{10}

\bibitem{armbruster_kinetic_2003}
D.~Armbruster, D.~Marthaler, and C.~Ringhofer.
\newblock Kinetic and {Fluid} {Model} {Hierarchies} for {Supply} {Chains}.
\newblock {\em Multiscale Modeling \& Simulation}, 2(1):43--61, January 2003.

\bibitem{barth_multilevel_2013}
Andrea Barth, Annika Lang, and Christoph Schwab.
\newblock Multilevel {Monte} {Carlo} method for parabolic stochastic partial differential equations.
\newblock {\em BIT Numerical Mathematics}, 53(1):3--27, March 2013.

\bibitem{bennett_uncertainty_2025}
W.~Bennett and R.~G. McClarren.
\newblock Uncertainty {Benchmarks} for {Time}-{Dependent} {Transport} {Problems}.
\newblock {\em Nuclear Science and Engineering}, 199(sup1):53, April 2025.

\bibitem{bennett_benchmarks_2022}
William Bennett and Ryan~G. McClarren.
\newblock Benchmarks for {Infinite} {Medium}, {Time} {Dependent} {Transport} {Problems} with {Isotropic} {Scattering}.
\newblock {\em Journal of Computational and Theoretical Transport}, 51(4):205--221, June 2022.

\bibitem{cao2015improved}
Wenhua Cao, Gino~J Lim, Yupeng Li, X~Ronald Zhu, and Xiaodong Zhang.
\newblock Improved beam angle arrangement in intensity modulated proton therapy treatment planning for localized prostate cancer.
\newblock {\em Cancers}, 7(2):574--584, 2015.

\bibitem{case_linear_1967}
Kenneth~M. Case and Paul~Frederick Zweifel.
\newblock {\em Linear {Transport} {Theory}}.
\newblock Addison-Wesley Publishing Company, 1967.

\bibitem{ceruti_robust_2024}
Gianluca Ceruti, Lukas Einkemmer, Jonas Kusch, and Christian Lubich.
\newblock A robust second-order low-rank {BUG} integrator based on the midpoint rule.
\newblock {\em BIT Numerical Mathematics}, 64(3):30, 2024.

\bibitem{ceruti_rank-adaptive_2022}
Gianluca Ceruti, Jonas Kusch, and Christian Lubich.
\newblock A rank-adaptive robust integrator for dynamical low-rank approximation.
\newblock {\em BIT Numerical Mathematics}, 62(4):1149--1174, December 2022.

\bibitem{ceruti_parallel_2024}
Gianluca Ceruti, Jonas Kusch, and Christian Lubich.
\newblock A {Parallel} {Rank}-{Adaptive} {Integrator} for {Dynamical} {Low}-{Rank} {Approximation}.
\newblock {\em SIAM Journal on Scientific Computing}, 46(3):B205--B228, June 2024.

\bibitem{ceruti_unconventional_2022}
Gianluca Ceruti and Christian Lubich.
\newblock An unconventional robust integrator for dynamical low-rank approximation.
\newblock {\em BIT Numerical Mathematics}, 62(1):23--44, March 2022.

\bibitem{chen_structure_2025}
Wei Chen, Giacomo Dimarco, and Lorenzo Pareschi.
\newblock Structure and asymptotic preserving deep neural surrogates for uncertainty quantification in multiscale kinetic equations, June 2025.
\newblock arXiv:2506.10636.

\bibitem{cliffe_multilevel_2011}
K~A Cliffe, M~B Giles, R~Scheichl, and A~L Teckentrup.
\newblock Multilevel {Monte} {Carlo} methods and applications to elliptic {PDEs} with random coefficients.
\newblock {\em Computing and Visualization in Science}, 14(1):3--15, August 2011.

\bibitem{collier_continuation_2015}
Nathan Collier, Abdul-Lateef Haji-Ali, Fabio Nobile, Erik Von~Schwerin, and Raúl Tempone.
\newblock A continuation multilevel {Monte} {Carlo} algorithm.
\newblock {\em BIT Numerical Mathematics}, 55(2):399--432, June 2015.

\bibitem{craft2014shared}
David Craft, Mark Bangert, Troy Long, D{\'a}vid Papp, and Jan Unkelbach.
\newblock Shared data for intensity modulated radiation therapy (imrt) optimization research: the cort dataset.
\newblock {\em GigaScience}, 3(1):2047--217X, 2014.

\bibitem{einkemmer_asymptotic-preserving_2024}
Lukas Einkemmer, Jingwei Hu, and Jonas Kusch.
\newblock Asymptotic-preserving and energy stable dynamical low-rank approximation.
\newblock {\em SIAM Journal on Numerical Analysis}, 62(1):73--92, 2024.

\bibitem{einkemmer_review_2025}
Lukas Einkemmer, Katharina Kormann, Jonas Kusch, Ryan~G. McClarren, and Jing-Mei Qiu.
\newblock A review of low-rank methods for time-dependent kinetic simulations.
\newblock {\em Journal of Computational Physics}, 538:114191, October 2025.

\bibitem{einkemmer_low-rank_2018}
Lukas Einkemmer and Christian Lubich.
\newblock A {Low}-{Rank} {Projector}-{Splitting} {Integrator} for the {Vlasov}--{Poisson} {Equation}.
\newblock {\em SIAM Journal on Scientific Computing}, 40(5):B1330--B1360, January 2018.

\bibitem{fairbanks_low-rank_2017}
Hillary~R. Fairbanks, Alireza Doostan, Christian Ketelsen, and Gianluca Iaccarino.
\newblock A low-rank control variate for multilevel {Monte} {Carlo} simulation of high-dimensional uncertain systems.
\newblock {\em Journal of Computational Physics}, 341:121--139, July 2017.

\bibitem{fu2023distributed}
Anqi Fu, Vicki~T Taasti, and Masoud Zarepisheh.
\newblock Distributed and scalable optimization for robust proton treatment planning.
\newblock {\em Medical physics}, 50(1):633--642, 2023.

\bibitem{ganapol_analytical_2008}
B~D. Ganapol and {OECD Nuclear Energy Agency.}
\newblock {\em Analytical benchmarks for nuclear engineering applications: case studies in neutron transport theory}.
\newblock Data bank. Nuclear Energy Agency, France, 2008.
\newblock OCLC: 430990895.

\bibitem{giles_multilevel_2008}
Michael~B. Giles.
\newblock Multilevel {Monte} {Carlo} {Path} {Simulation}.
\newblock {\em Operations Research}, 56(3):607--617, June 2008.

\bibitem{giles_multilevel_2015}
Michael~B. Giles.
\newblock Multilevel {Monte} {Carlo} methods.
\newblock {\em Acta Numerica}, 24:259--328, May 2015.

\bibitem{grote_uncertainty_2022}
Marcus~J. Grote, Simon Michel, and Fabio Nobile.
\newblock Uncertainty {Quantification} by {Multilevel} {Monte} {Carlo} and {Local} {Time}-{Stepping} for {Wave} {Propagation}.
\newblock {\em SIAM/ASA Journal on Uncertainty Quantification}, 10(4):1601--1628, December 2022.

\bibitem{haji-ali_optimization_2016}
Abdul-Lateef Haji-Ali, Fabio Nobile, Erik von Schwerin, and Raúl Tempone.
\newblock Optimization of mesh hierarchies in multilevel {Monte} {Carlo} samplers.
\newblock {\em Stochastics and Partial Differential Equations Analysis and Computations}, 4(1):76--112, March 2016.

\bibitem{jin_asymptotic-preserving_2024}
Shi Jin, Zheng Ma, and Keke Wu.
\newblock Asymptotic-{Preserving} {Neural} {Networks} for {Multiscale} {Kinetic} {Equations}.
\newblock {\em Communications in Computational Physics}, 35(3):693--723, January 2024.

\bibitem{kieri_discretized_2016}
Emil Kieri, Christian Lubich, and Hanna Walach.
\newblock Discretized {Dynamical} {Low}-{Rank} {Approximation} in the {Presence} of {Small} {Singular} {Values}.
\newblock {\em SIAM Journal on Numerical Analysis}, 54(2):1020--1038, January 2016.

\bibitem{koch_dynamical_2007}
Othmar Koch and Christian Lubich.
\newblock Dynamical {Low}‐{Rank} {Approximation}.
\newblock {\em SIAM Journal on Matrix Analysis and Applications}, 29(2):434--454, January 2007.

\bibitem{koellermeier_macro-micro_2024}
Julian Koellermeier, Philipp Krah, and Jonas Kusch.
\newblock Macro-micro decomposition for consistent and conservative model order reduction of hyperbolic shallow water moment equations: a study using {POD}-{Galerkin} and dynamical low-rank approximation.
\newblock {\em Advances in Computational Mathematics}, 50(4):76, August 2024.

\bibitem{koellermeier2020analysis}
Julian Koellermeier and Marvin Rominger.
\newblock Analysis and numerical simulation of hyperbolic shallow water moment equations.
\newblock {\em Communications in Computational Physics}, 28(3):1038--1084, 2020.

\bibitem{kowalski_moment_2019}
Julia Kowalski and Manuel Torrilhon.
\newblock Moment {Approximations} and {Model} {Cascades} for {Shallow} {Flow}.
\newblock {\em Communications in Computational Physics}, 25(3):669--702, January 2019.

\bibitem{kusch_second-order_2025}
Jonas Kusch.
\newblock Second-order robust parallel integrators for dynamical low-rank approximation.
\newblock {\em BIT Numerical Mathematics}, 65(3):31, September 2025.

\bibitem{kusch_dynamical_2022}
Jonas Kusch, Gianluca Ceruti, Lukas Einkemmer, and Martin Frank.
\newblock Dynamical {Low}-{Rank} {Approximation} for {Burgers}' {Equation} with {Uncertainty}.
\newblock {\em International Journal for Uncertainty Quantification}, 12(5):1--21, 2022.

\bibitem{kusch_stability_2023}
Jonas Kusch, Lukas Einkemmer, and Gianluca Ceruti.
\newblock On the {Stability} of {Robust} {Dynamical} {Low}-{Rank} {Approximations} for {Hyperbolic} {Problems}.
\newblock {\em SIAM Journal on Scientific Computing}, 45(1):A1--A24, February 2023.

\bibitem{kusch_kit-rt_2023}
Jonas Kusch, Steffen Schotthöfer, Pia Stammer, Jannick Wolters, and Tianbai Xiao.
\newblock {KiT}-{RT}: {An} {Extendable} {Framework} for {Radiative} {Transfer} and {Therapy}.
\newblock {\em ACM Trans. Math. Softw.}, 49(4):38:1--38:24, December 2023.

\bibitem{kusch_robust_2023}
Jonas Kusch and Pia Stammer.
\newblock A robust collision source method for rank adaptive dynamical low-rank approximation in radiation therapy.
\newblock {\em ESAIM: Mathematical Modelling and Numerical Analysis}, 57(2):865--891, March 2023.

\bibitem{leveque_finite_2002}
Randall~J. LeVeque.
\newblock {\em Finite {Volume} {Methods} for {Hyperbolic} {Problems}}.
\newblock Cambridge University Press, 1 edition, August 2002.

\bibitem{lewis_computational_1984}
Elmer~Eugene Lewis and Warren~F Miller.
\newblock {\em Computational methods of neutron transport}.
\newblock United States: John Wiley and Sons, Inc., 1984.

\bibitem{liu_bi-fidelity_2020}
Liu Liu and Xueyu Zhu.
\newblock A bi-fidelity method for the multiscale {Boltzmann} equation with random parameters.
\newblock {\em Journal of Computational Physics}, 402:108914, February 2020.

\bibitem{lomax2008intensitya}
AJ~Lomax.
\newblock Intensity modulated proton therapy and its sensitivity to treatment uncertainties 1: the potential effects of calculational uncertainties.
\newblock {\em Physics in Medicine \& Biology}, 53(4):1027--1042, 2008.

\bibitem{lomax2008intensityb}
AJ~Lomax.
\newblock Intensity modulated proton therapy and its sensitivity to treatment uncertainties 2: the potential effects of inter-fraction and inter-field motions.
\newblock {\em Physics in Medicine \& Biology}, 53(4):1043--1056, 2008.

\bibitem{lubich_projector-splitting_2014}
Christian Lubich and Ivan~V. Oseledets.
\newblock A projector-splitting integrator for dynamical low-rank approximation.
\newblock {\em BIT Numerical Mathematics}, 54(1):171--188, March 2014.

\bibitem{luo_multilevel_2019}
Yan Luo and Zhu Wang.
\newblock A {Multilevel} {Monte} {Carlo} {Ensemble} {Scheme} for {Random} {Parabolic} {PDEs}.
\newblock {\em SIAM Journal on Scientific Computing}, 41(1):A622--A642, January 2019.

\bibitem{mishra_sparse_2012}
S.~Mishra and Ch. Schwab.
\newblock Sparse tensor multi-level {Monte} {Carlo} finite volume methods for hyperbolic conservation laws with random initial data.
\newblock {\em Mathematics of Computation}, 81(280):1979--2018, April 2012.

\bibitem{mishra_multi-level_2012}
S.~Mishra, Ch. Schwab, and J.~Šukys.
\newblock Multi-level {Monte} {Carlo} finite volume methods for nonlinear systems of conservation laws in multi-dimensions.
\newblock {\em Journal of Computational Physics}, 231(8):3365--3388, April 2012.

\bibitem{musharbash_dual_2018}
Eleonora Musharbash and Fabio Nobile.
\newblock Dual {Dynamically} {Orthogonal} approximation of incompressible {Navier} {Stokes} equations with random boundary conditions.
\newblock {\em Journal of Computational Physics}, 354:135--162, February 2018.

\bibitem{musharbash_symplectic_2020}
Eleonora Musharbash, Fabio Nobile, and Eva Vidličková.
\newblock Symplectic dynamical low rank approximation of wave equations with random parameters.
\newblock {\em BIT Numerical Mathematics}, 60(4):1153--1201, 2020.

\bibitem{nobile_high-order_2026}
Fabio Nobile and Sébastien Riffaud.
\newblock High-{Order} {BUG} {Dynamical} {Low}-{Rank} {Integrators} {Based} on {Explicit} {Runge}–{Kutta} {Methods}.
\newblock {\em Journal of Scientific Computing}, 107(3):102, June 2026.

\bibitem{pareschi_introduction_2021}
Lorenzo Pareschi.
\newblock An {Introduction} to {Uncertainty} {Quantification} for {Kinetic} {Equations} and {Related} {Problems}.
\newblock In Giacomo Albi, Sara Merino-Aceituno, Alessia Nota, and Mattia Zanella, editors, {\em Trails in {Kinetic} {Theory}: {Foundational} {Aspects} and {Numerical} {Methods}}, pages 141--181. Springer International Publishing, Cham, 2021.

\bibitem{patwardhan_parallel_2026}
Chinmay Patwardhan and Jonas Kusch.
\newblock A {Parallel}, {Energy}-{Stable} {Low}-{Rank} {Integrator} for {Nonlinear} {Multi}-{Scale} {Thermal} {Radiative} {Transfer}.
\newblock {\em Journal of Computational and Theoretical Transport}, 55(2):145--189, February 2026.

\bibitem{patwardhan_low-rank_2026}
Chinmay Patwardhan, Pia Stammer, Emil Løvbak, Jonas Kusch, and Sebastian Krumscheid.
\newblock Low-{Rank} {Variance} {Reduction} for {Uncertain} {Radiative} {Transfer} with {Control} {Variates}.
\newblock In Christiane Lemieux and Ben Feng, editors, {\em Monte {Carlo} and {Quasi}-{Monte} {Carlo} 2024}, volume 522, pages 357--375. Springer Nature Switzerland, Cham, 2026.
\newblock Series Title: Springer Proceedings in Mathematics \& Statistics.

\bibitem{perko2016fast}
Zolt{\'a}n Perk{\'o}, Sebastian~R Van Der~Voort, Steven Van De~Water, Charlotte~MH Hartman, Mischa Hoogeman, and Danny Lathouwers.
\newblock Fast and accurate sensitivity analysis of impt treatment plans using polynomial chaos expansion.
\newblock {\em Physics in Medicine \& Biology}, 61(12):4646--4664, 2016.

\bibitem{prugger_dynamical_2023}
Martina Prugger, Lukas Einkemmer, and Carlos~F. Lopez.
\newblock A dynamical low-rank approach to solve the chemical master equation for biological reaction networks.
\newblock {\em Journal of Computational Physics}, 489:112250, September 2023.

\bibitem{sapsis_dynamically_2009}
Themistoklis~P Sapsis and Pierre~FJ Lermusiaux.
\newblock Dynamically orthogonal field equations for continuous stochastic dynamical systems.
\newblock {\em Physica D: Nonlinear Phenomena}, 238(23-24):2347--2360, 2009.

\bibitem{schotthofer_reference_2025}
Steffen Schotthöfer and Cory Hauck.
\newblock Reference solutions for linear radiation transport: the {Hohlraum} and {Lattice} benchmarks, May 2025.

\bibitem{stammer2026high}
Pia Stammer, Niklas Wahl, Jonas Kusch, and Danny Lathouwers.
\newblock A high-order deterministic dynamical low-rank method for proton transport in heterogeneous media.
\newblock {\em Journal of Computational Physics}, page 114879, 2026.

\bibitem{ueckermann2013numerical}
M.P. Ueckermann, P.F.J. Lermusiaux, and T.P. Sapsis.
\newblock Numerical schemes for dynamically orthogonal equations of stochastic fluid and ocean flows.
\newblock {\em Journal of Computational Physics}, 233:272--294, 2013.

\bibitem{welford_note_1962}
B.~P. Welford.
\newblock Note on a {Method} for {Calculating} {Corrected} {Sums} of {Squares} and {Products}.
\newblock {\em Technometrics}, 4(3):419--420, August 1962.

\end{thebibliography}

\end{document}